\documentclass[12pt, reqno]{amsart}
\usepackage{amsmath, amsthm, amscd, amsfonts, amssymb, graphicx,color,tikz}
\usepackage[bookmarksnumbered, colorlinks, plainpages]{hyperref}
\usepackage{mathrsfs}
\usepackage{enumerate}
\usepackage[utf8]{inputenc}
\usepackage{bbold}
\usepackage{booktabs}
\usepackage{tikz}

\usepackage[all]{xy}  

\newtheorem{theorem}{Theorem}[section]
\newtheorem{lemma}[theorem]{Lemma}
\newtheorem{proposition}[theorem]{Proposition}
\newtheorem{corollary}[theorem]{Corollary}
\theoremstyle{definition}

\newtheorem{example}[theorem]{Example}

\numberwithin{equation}{section}

\newcommand{\s}{\mathrm{span\,}}

\newcommand{\CC}{\mathbb C}

\newcommand{\DD}{\mathbb D}
\newcommand{\TT}{\mathbb T}
\newcommand{\ZZ}{\mathbb Z}
\newcommand{\NN}{\mathbb N}
\newcommand{\diag}{\mathrm{diag\,}}

\newcommand{\alg}{\mathrm{alg\,}}

\newcommand{\dist}{\mathrm{dist\,}}
\newcommand{\supp}{\mathrm{supp\,}}
\newcommand{\PP}{\mathcal P}
\newcommand{\Tr}{\mathrm{Tr\,}}
\newcommand{\Ee}{\mathcal E}
\newcommand{\Hh}{\mathcal H}
\newcommand{\Mcal}{\mathcal M}

\newcommand{\Rat}{\mathrm{Rat\,}}
\newcommand{\UU}{\mathbb U}

\newcommand{\fcd}{\mathcal{F}({\mathbb{C}}^d)}

\allowdisplaybreaks[4]

\hypersetup{colorlinks=true,linkcolor=blue,citecolor=green}

\begin{document}
\setcounter{page}{1}

\title[Commutants of diagonal operators]{Commutants of diagonal operators}
\date{\today}

\author[M. Wang and X. Yao]{Maofa Wang and Xingxing Yao}

\address{School of Mathematics and Statistics, Wuhan University, Wuhan 430072, China.} \email{mfwang.math@whu.edu.cn}

\address{School of Mathematics and Physics, Wuhan Institute of Technology, Wuhan 430205, China.}
\email{xxyao.math@wit.edu.cn}




\begin{abstract}
The aim of this paper is devoted to the profound and comprehensive characterization when diagonal operators having the minimal, tiny or double commutant property on a separable Hilbert space, in terms of the topological and geometric properties of their spectra.
\end{abstract}
\maketitle


\section{Introduction}

Denote by $\mathcal{B}(H)$ the algebra of all bounded linear operators on a Hilbert space $H$.
Given any subset $\mathcal{A}\subseteq \mathcal{B}(H)$,
recall that the commutant of $\mathcal{A}$ is defined as the family of all operators that commute with all the operators in $\mathcal{A}$,
that is,
\begin{equation*}
  \mathcal{A}':=\bigcap_{A\in\mathcal{A}}\{A\}'=\bigcap_{A\in\mathcal{A}}\{T\in\mathcal{B}(H):TA=AT\},
\end{equation*}
a closed subalgebra of $\mathcal{B}(H)$ in the weak operator topology.
The double commutant of $\mathcal{A}$ is defined as $\mathcal{A}''=(\mathcal{A}')'$.
Clearly, $\mathcal{A}$ is contained in $\mathcal{A}''$,
and then so is $\overline{\mathcal{A}}^{WOT}$,
the closure of $\mathcal{A}$ in the weak operator topology.
As is well-known,
$\overline{\mathcal{A}}^{WOT}=\mathcal{A}''$ holds for any self-adjoint, unital subalgebra $\mathcal{A}\subseteq \mathcal{B}(H)$,
by the von Neumann's double commutant theorem. A subalgebra $\mathcal{A}\subseteq \mathcal{B}(H)$ is said to possess the \emph{double commutant property} provided that $\overline{\mathcal{A}}^{WOT}=\mathcal{A}''$. Similarly, we say that a subalgebra $\mathcal{A}\subseteq \mathcal{B}(H)$ has the \emph{minimal commutant property} provided that $\overline{\mathcal{A}}^{WOT}=\{A\}'$. It is clear that a subalgebra with a minimal commutant has the double commutant property.

Denote by $\PP$ the set of all polynomials. Given a single operator $A\in \mathcal{B}(H)$,
let $\alg (A)$ stand for the unital algebra generated by the operator $A$, that is,
\begin{equation*}
  \alg (A):=\{p(A):p\in \PP\}.
\end{equation*}
Moreover, denote by $\alg^{a} (A)$ the set of $f(A)$ with rational functions $f$ whose poles are outside the spectrum of $A$. It is easy to see that both $\overline{\alg (A)}^{WOT}$ and $\overline{\alg^{a} (A)}^{WOT}$ are commutative algebras such that
\begin{equation*}
  \overline{\alg (A)}^{WOT}\subseteq \overline{\alg^{a} (A)}^{WOT}\subseteq \{A\}''\subseteq\{A\}'.
\end{equation*}
One can define $A$ has the \emph{minimal commutant property} or possess the \emph{double commutant property}, too. Thanks to Herrero's insightful work \cite{He}, $\big\{A\in\mathcal{B}(H):\overline{\alg^{a} (A)}^{WOT}=\{A\}'\big\}$ is norm-dense in $\mathcal{B}(H)$. Similarly, we say that $A$ has the \emph{tiny commutant property} (resp. the \emph{tiny double commutant property}), if $\overline{\alg^{a} (A)}^{WOT}=\{A\}'$ (resp. $\overline{\alg^{a} (A)}^{WOT}=\{A\}''$). That has already been investigated in several special operators, for example, the matrix on the finite-dimensional spaces (a matrix has a minimal commutant if and only if it is nonderogatory), a forward weighted shift on $\ell^2$ or the discrete Ces\`{a}ro operator (see \cite{SW}), and more recently, composition operators on various spaces (see \cite{BWY,BY,LLSR18,LLSR19}).

Now this paper seeks to characterize thoroughly when a diagonal operator having the minimal or double commutant property on a separable Hilbert space. As is well-known, the minimal commutant property of a direct sum implies that of each direct summand; but the converse fails in general. Fortunately, we can provide a sufficient criterion when the minimal commutant property of a direct sum can be deduced from that of all direct summands in Section \ref{sec:com-d-s} (see Theorem \ref{thm-sum-com}).
Section \ref{sec:com-diag} intends to characterize when a diagonal operator has the minimal or double commutant property on a separable Hilbert space by using the topological property of its spectrum (see Theorem \ref{thm-diag} and Theorem \ref{thm-diag-DCP}). Some interesting particular cases, including Corollary \ref{cor-diag-real}, Example \ref{ex-diag-ce} and Corollary \ref{cor:J-one}, are also showed. On the other hand, some easily verifiable criteria are given via the geometric properties of the spectrum of the diagonal operator, by using the complex-polynomial approximation (see Theorem \ref{thm:sepc-mc}). The similar characterization for the tiny or double commutant property is given in Section \ref{sec:tiny-com-diag} (see Theorem \ref{thm-diag-tiny} and Theorem \ref{thm-diag-tiny-double}). As an application, the corresponding properties of some composition operators are considered in Section \ref{sec:app-co} (see Theorem \ref{ex-diag-mc},  Corollaries \ref{MC-com-Fock} and \ref{MC-com-J-Fock}).

\section{Commutant of direct sums}\label{sec:com-d-s}

We first collect some properties about the relation between the commutant of a direct sum and that of direct summands from \cite{LLSR18,LLSR19}.

\begin{lemma}\label{sum-com}
Let $H=H_{1}\oplus H_{2}$, let $A_{j}\in \mathcal{B}(H_{j})$ for $j=1,2$, and let $A=A_{1}\oplus A_{2}$.
\begin{itemize}
  \item[(1)] If $A$ has a minimal commutant, then so do both $A_{1}$ and $A_{2}$.
  \item[(2)] If $\|A_{2}\|<1$, then $\overline{\alg(I_{1}\bigoplus A_{2})}^{\sigma}=\alg(I_{1})\bigoplus\overline{\alg(A_{2})}^{\sigma}$.
\end{itemize}
Moreover, for subalgebras $\mathcal{A}_{j}\subseteq\mathcal{B}(H_{j})$ for $j=1$ or $2$, $(\mathcal{A}_{1}\bigoplus\mathcal{A}_{2})'=\mathcal{A}_{1}'\bigoplus\mathcal{A}_{2}'$ and $(\mathcal{A}_{1}\bigoplus\mathcal{A}_{2})''=\mathcal{A}_{1}''\bigoplus\mathcal{A}_{2}''$.
\end{lemma}

Conversely, we have

\begin{lemma}\label{sum-com-inv}
Let $H=H_{1}\oplus H_{2}$, let $A_{j}\in \mathcal{B}(H_{j})$ for $j=1,2$, and let $A=A_{1}\oplus A_{2}$ with $\sigma(A_1)\cap\sigma(A_2)=\emptyset$. Then
$\{A_1\oplus A_2\}'=\{A_1\}'\oplus\{A_2\}'.$
\end{lemma}

\begin{proof}
For any $X\in \mathcal{B}(H)$, write its block decomposition
\begin{equation*}
X=\begin{pmatrix}
X_{11}&X_{12}\\
X_{21}&X_{22}
\end{pmatrix}.
\end{equation*}
The condition $AX=XA$ is equivalent to the system
\begin{equation*}
\begin{cases}
A_1X_{11}=X_{11}A_1,\\
A_1X_{12}=X_{12}A_2,\\
A_2X_{21}=X_{21}A_1,\\
A_2X_{22}=X_{22}A_2.
\end{cases}
\end{equation*}
To be proceed, recall the Rosenblum operator
\begin{align*}
\tau_{A_1,A_2}:\mathcal{B}(H_2,H_1) & \longrightarrow \mathcal{B}(H_2,H_1),  \\
  X & \longmapsto A_1X-XA_2
\end{align*}
and Sylvester-Rosenblum Theorem states that
$$\sigma(A_1)\cap\sigma(A_2)=\emptyset \
\Longrightarrow \ \tau_{A_1,A_2}\text{ is invertible},$$
(for example see \cite{BR}). Thus, if $\sigma(A_1)\cap\sigma(A_2)=\varnothing,$
then $A_1X_{12}=X_{12}A_2$ and $A_2X_{21}=X_{21}A_1$ implies that $X_{12}=0$ and $X_{21}=0$, respectively. Thus $\{A_1\oplus A_2\}'=\{A_1\}'\oplus\{A_2\}'.$
\end{proof}

Although‌ $\{A_1\oplus A_2\}'=\{A_1\}'\oplus\{A_2\}'$ by the property $\sigma(A_1)\cap\sigma(A_2)=\emptyset$, the minimal commutant property of $A=A_1\oplus A_2$  cannot be directly deduced from that of $A_1$ and $A_2$, which may be attributed to the fact that we require a single common polynomial to approximate both blocks simultaneously and separate polynomials for each block are insufficient. We now present a general sufficient condition.

\begin{theorem}\label{thm-sum-com}
Let $H=H_{1}\oplus H_{2}$, let $A_{j}\in \mathcal{B}(H_{j})$ for $j=1,2$, and let $A=A_{1}\oplus A_{2}$ with $\sigma(A_1)\cap\sigma(A_2)=\emptyset$. If both $A_1\in \mathcal{B}(H_1)$ and $A_2\in \mathcal{B}(H_2)$ possess the minimal commutant property, and there exists a sequence of polynomials $r_n$ such that
$$\|r_n(A_1)-I_{H_1}\|\longrightarrow0,
\qquad\|r_n(A_2)\|\longrightarrow0.$$
Then $A$ has a minimal commutant.
\end{theorem}

\begin{proof}
Since $\sigma(A_1)\cap\sigma(A_2)=\emptyset$, then $\{A_1\oplus A_2\}'=\{A_1\}'\oplus\{A_2\}'$ by Lemma \ref{sum-com-inv}.
Given arbitrary $T=T_1\oplus T_2 \in\{A_1\oplus A_2\}'.$ Since both $A_1$ and $A_2$ have minimal commutants, there are polynomials $p_k,q_m$ such that
\begin{equation*}
p_k(A_1)\rightarrow T_1 \mbox{ and }
q_m(A_2)\rightarrow T_2
\end{equation*}
in the weak operator topology.

Note that $p_k(A_1)\oplus q_m(A_2)$ is not itself a polynomial in $A_1\oplus A_2$. We construct a unified polynomial by setting $P_{n,k,m}(z)=r_n(z)p_k(z)+\big(1-r_n(z)\big)q_m(z).$ Then
$$\begin{aligned}
P_{n,k,m}(A_1)
&=r_n(A_1)p_k(A_1)
 +\big(I_{H_1}-r_n(A_1)\big)q_m(A_1),\\
P_{n,k,m}(A_2)
&=r_n(A_2)p_k(A_2)
 +\big(I_{H_2}-r_n(A_2)\big)q_m(A_2).
\end{aligned}$$
Fix $k,m$, and take $n$ sufficiently large. By the hypotheses,
$$P_{n,k,m}(A_1)\longrightarrow p_k(A_1),\qquad P_{n,k,m}(A_2)\longrightarrow q_m(A_2)$$
in the operator norm topology. Furthermore, for any $x,y\in H$, writing $x=x_1\oplus x_2$ and $y=y_1\oplus y_2$, we have
\begin{equation*}
 \langle P_{n,k,m}(A)x,y\rangle=\langle P_{n,k,m}(A_1)x_1,y_1\rangle_{H_1}+\langle P_{n,k,m}(A_2)x_2,y_2\rangle_{H_2}.
\end{equation*}
Then for sufficiently large $k,m$, $P_{n,k,m}(A)$ converges weakly to $T=T_1\oplus T_2$. Thus $A$ has a minimal commutant.
\end{proof}

The above polynomial separation condition for the two operators is not always easy to check directly. We give a usefully sufficient criterion for normal operators.

\begin{corollary}\label{sum-normal}
Given $m\in\NN$, let $A_j\in \mathcal{B}(H_j)$ be normal operators having a minimal commutant for $1\leq j\leq m$. If $\sigma(A_{j})\cap \sigma(A_k)=\emptyset$ for any $j\neq k$, and $\CC\setminus \cup_{j=1}^{m}\sigma(A_j)$ is connected, then $\oplus_{j=1}^{m}A_j$ has a minimal commutant.
\end{corollary}

\begin{proof}
Since the direct sum of two normal operators is also normal and the spectrum of the direct sum is equal to the union set of spectrum of summand operators, then it suffices to consider the two-term direct sum case. For this, denote $K_1=\sigma(A_1),$ $K_2=\sigma(A_2)$ and $K=K_1\cup K_2$. Then $K_1,K_2$ are disjoint compact sets, and thus we may define a smooth function on an open neighbourhood of $K$ by
\begin{equation*}
f(z)=
\begin{cases}
1,&z\ \text{in a neighbourhood of }K_1,\\
0,&z\ \text{in a neighbourhood of }K_2.
\end{cases}
\end{equation*}
By Runge’s theorem (for example see \cite{Ga}), there exists a sequence of polynomials $r_n$ such that
\begin{equation*}
\sup_{z\in K_1}|r_n(z)-1|\to0 \mbox{ and }
\sup_{z\in K_2}|r_n(z)|\to0.
\end{equation*}
Because $A_1,A_2$ are normal, the continuous functional calculus for normal operators gives
\begin{equation*}
\|r_n(A_1)-I\|
=\sup_{z\in K_1}|r_n(z)-1|\to0,
\end{equation*}
and
\begin{equation*}
\|r_n(A_2)\|
=\sup_{z\in K_2}|r_n(z)|\to0,
\end{equation*}
(for example see \cite[Chapter 5]{Mac}). Then the desired result can be immediately obtained following from Theorem \ref{thm-sum-com}.
\end{proof}

\noindent \textbf{Note}.  In general, the convergence $\sup_{z\in\sigma(A_j)}|r_n(z)-c|\longrightarrow0$ does not imply $\|r_n(A_j)-cI\|\longrightarrow0.$ Accordingly, the normality assumption in the corollary cannot be omitted.

\begin{example}
Let $H_1=H_2=\mathbb C^2,$ and set
\begin{equation*}
A_1=A_2=
\begin{pmatrix}
0&0\\
0&1
\end{pmatrix}.
\end{equation*}
Then both $A_1$ and $A_2$ have minimal commutants, but $A_1\oplus A_2$ fails to have a minimal commutant.
\end{example}

Indeed, it follows immediately from the fact that if $H$ is finite-dimensional, $A\in \mathcal{B}(H)$ has a minimal commutant if and only if it is nonderogatory, i.e., its eigenvalues have geometric multiplicity equal to $1$. Moreover, in terms of the Rosenblum operator, $\tau_{A_1,A_2}(Y)=A_1Y-YA_2=0$ holds for all $Y\in \mathcal{M}_{2}(\CC)$. And the counterexample shows that the existence of non-zero off-diagonal commuting operators is exactly what destroys the minimal commutant property.

\section{Commutant of diagonal operators}\label{sec:com-diag}

\subsection{Criteria}~
\medskip

We fist give criteria when a diagonal operator has the minimal or double commutant property on a separable Hilbert space by using the topological property of its spectrum. For this, we recall the following result about annihilators from \cite[Chapter 4, Theorem 4.7(b)]{Rudin91}.

\begin{lemma}\label{anni-dual}
Let $X$ be a Banach space and $X^\ast$ be its dual space. For $M\subseteq X^\ast$, let the annihilator space of $M$ be ${}_\perp M=\{x\in X:f(x)=0, \mbox{ for all } f\in M\}.$
Then $\overline M^{\,\sigma(X^*,X)}=({}_\perp M)^\perp,$
where
\begin{equation*}
({}_\perp M)^\perp=\{f\in X^*:f(x)=0, \mbox{ for all } x\in{}_\perp M\}.
\end{equation*}
Especially,
\begin{equation*}
\overline M^{\,\sigma(X^*,X)}=X^*\iff {}_\perp M=\{0\}.
\end{equation*}
\end{lemma}

\begin{theorem}\label{thm-diag}
Let $D=\diag(d_n)$ be a bounded diagonal operator on a separable Hilbert space $H$. Then $D$ has the minimal commutant property if and only if the sequence $(d_n)$ is pairwise distinct, and the only sequence $c=(c_n)\in\ell^1$ satisfying
\begin{equation*}
\sum_{n=1}^{\infty}c_nd_n^k=0, \quad k=0,1,2,\ldots,
\end{equation*}
is the zero sequence.

In other words, the minimal commutant property of a diagonal operator with pairwise distinct diagonal entries is equivalent to uniqueness of the corresponding \(\ell^1\)-moment problem.
\end{theorem}

\begin{proof}
If $D$ has the minimal commutant property, then $(d_n)$ is pairwise distinct, and thus
\begin{equation*}
\{D\}'=\left\{\diag(\mu_n):(\mu_n)\in\ell^\infty\right\}\subseteq \mathcal{B}(H).
\end{equation*}
Consider the isometric isomorphism
\begin{equation*}
\Phi:\ell^\infty\longrightarrow \mathcal{B}(H),\quad
\Phi(\mu)=\operatorname{diag}(\mu_n).
\end{equation*}
Then $\Phi(\ell^\infty)=\{D\}'.$ We claim that the weak operator topology on $\Phi(\ell^\infty)$ is equivalent to $\sigma(\ell^{\infty},\ell^{1})$, the weak $\ast$ topology of $\ell^{\infty}$ defined by a net $(\mu^{(\alpha)})\subseteq \ell^{\infty}$ converging to $\mu\in \ell^{\infty}$ equivalent to
\begin{equation*}
  \sum_{n=1}^{\infty}c_n\mu_n^{(\alpha)}\longrightarrow\sum_{n=1}^{\infty}c_n\mu_n\ \mbox{ for any } c=(c_n)\in\ell^1.
\end{equation*}
Indeed, given $f=\sum_{n=1}^{\infty} \xi_ne_n, g=\sum_{n=1}^{\infty}\eta_ne_n\in H$, define $c_n=\xi_n\overline{\eta_n}$ and then $c=(c_n)\in\ell^{1}$ by the Cauchy-Schwarz inequality. Particularly, for any $T=\diag(\mu_n)\in \Phi(\ell^\infty)$ with $(\mu_n)\in\ell^\infty$,
\begin{equation*}
 \langle Tf,g\rangle = \sum_{n=1}^{\infty}\mu_nc_n.
\end{equation*}
Also, any linear functional on $\mathcal B(H)$ which is continuous with respect to the weak operator topology can be expressed as a linear combination of finitely many linear functionals of the form $\langle Tf,g\rangle$ with $f,g\in H$ (see \cite[Chapter 2]{Ta}). So the weak convergence on $\Phi(\ell^\infty)$ means that a net $(\mu^{(\alpha)})\subseteq \ell^{\infty}$ converges to $\mu\in\ell^{\infty}$, i.e.,
\begin{equation*}
  \sum_{n=1}^{\infty}c_n\mu_n^{(\alpha)}\longrightarrow\sum_{n=1}^{\infty}c_n\mu_n\ \mbox{ for any } c=(c_n)\in\ell^1.
\end{equation*}
Conversely, given any $c=(c_n)\in\ell^1$, set $\xi_n=\sqrt{|c_n|}$ and
\begin{equation*}
\eta_n=\begin{cases}\dfrac{\overline{c_n}}{\sqrt{|c_n|}},&c_n\ne0,\\[2mm]
0,&c_n=0.
\end{cases}
\end{equation*}
Then
\begin{equation*}
\sum_{n=1}^{\infty}|\xi_n|^2=\sum_{n=1}^{\infty}|c_n|<\infty, \mbox{ and }
\sum_{n=1}^{\infty}|\eta_n|^2=\sum_{n=1}^{\infty}|c_n|<\infty.
\end{equation*}
So $\xi,\eta\in\ell^2$ and $\xi_n\overline{\eta_n}=c_n,$ which in turn implies that
\begin{equation*}
\langle \Phi(\mu)f,g\rangle=\sum_{n=1}^{\infty}c_n\mu_n\ \ \mbox{ with } f=\sum_{n=1}^{\infty} \xi_ne_n, g=\sum_{n=1}^{\infty}\eta_ne_n\in H.
\end{equation*}

Due to $\sup_n|d_n|=\|D\|<\infty$, for any $p\in \PP$, we have $\sup_n|p(d_n)|<\infty$, and then
\begin{equation*}
\ell^{\infty}_D=\left\{(p(d_n))_{n\geq1}:p\in\PP\right\}\subseteq\ell^{\infty}.
\end{equation*}
Clearly, $\Phi(\ell^{\infty}_D)=\alg\{D\}.$ Therefore, $D$ has a minimal commutant if and only if
\begin{equation}\label{equi-relation}
  \overline{\Phi(\ell^{\infty}_D)}^{WOT}=\overline{\alg\{D\}}^{WOT}=\{D\}'=\Phi(\ell^{\infty}),
\end{equation}
which is equivalent to
\begin{equation}\label{ell-relation}
\overline{\ell^{\infty}_D}^{\sigma(\ell^\infty,\ell^1)}
=\ell^\infty.
\end{equation}
According to the dual relation $(\ell^{1})^{\ast}=\ell^{\infty}$, by Lemma \ref{anni-dual},
\begin{equation}\label{ell-relation-anni}
\overline{\ell^{\infty}_D}^{\,\sigma(\ell^\infty,\ell^1)}=\ell^\infty\
\mbox{ if and only if }\
{}_\perp\ell^{\infty}_D=\{0\}.
\end{equation}

To be proceed, note that the annihilator space of $\ell^{\infty}_D$ is
\begin{equation*}
{}_\perp\ell^{\infty}_D=
\left\{c\in\ell^1:\sum_{n=1}^{\infty}c_n p(d_n)=0 \mbox{ for all } p\in \PP\right\}.
\end{equation*}
It is easy to see that
\begin{equation*}
{}_\perp\ell^{\infty}_D=
\left\{c\in\ell^1:\sum_{n=1}^{\infty}c_nd_n^k=0  \mbox{ for all } k\geq0\right\},
\end{equation*}
due to the fact $(d_n)\in\ell^{\infty}$. That, together with the relations \eqref{equi-relation}-\eqref{ell-relation-anni}, implies that
\begin{equation*}
  \overline{\alg\{D\}}^{WOT}=\{D\}'\Longleftrightarrow \left\{c\in\ell^1:\sum_{n=1}^{\infty}c_nd_n^k=0  \mbox{ for all } k\geq0\right\}=\{0\},
\end{equation*}
as we wanted.
\end{proof}

Similar to the minimal commutant property in Theorem \ref{thm-diag}, the following holds for the double commutant property.

\begin{theorem}\label{thm-diag-DCP}
Let $D=\diag(d_n)$ be a bounded diagonal operator on a separable Hilbert space $H$, and let $\Lambda=\{d_n:n\geq 1\}$ be the set of its distinct diagonal entries. Then $D$ has the double commutant property if and only if
\begin{equation*}
  \left\{c\in\ell^1(\Lambda):\sum_{\lambda\in\Lambda}c_\lambda \lambda^k=0\text{ for every }k\geq0\right\}=\{0\}.
\end{equation*}
\end{theorem}

\begin{proof}
For each \(\lambda\in\Lambda\), define $I_\lambda=\{n\geq 1:d_n=\lambda\}$ and let $H_\lambda=\overline{\operatorname{span}}\{e_n:n\in I_\lambda\}.$
Then
\begin{equation*}
H=\bigoplus_{\lambda\in\Lambda}H_\lambda, \mbox{ and }
D=\bigoplus_{\lambda\in\Lambda}\lambda I_{H_\lambda}.
\end{equation*}
By Lemma \ref{sum-com-inv},
\begin{equation*}
\{D\}'=\left\{\bigoplus_{\lambda\in\Lambda}T_\lambda:
T_\lambda\in \mathcal{B}(H_\lambda),\ \sup_{\lambda\in\Lambda}\|T_\lambda\|<\infty\right\}.
\end{equation*}
Let $P_\lambda$ denote the orthogonal projection from $H$ onto $H_\lambda$. Then each $P_\lambda$ belongs to $\{D\}'$, and thus every operator $S\in\{D\}''$ commutes with all $P_\lambda$. So $S$ is block diagonal and can be written as
\begin{equation*}
S=\bigoplus_{\lambda\in\Lambda}S_\lambda,\quad S_\lambda\in \mathcal{B}(H_\lambda).
\end{equation*}

Moreover, for any fixed \(\lambda\in\Lambda\), the commutant \(\{D\}'\) contains every operator acting arbitrarily on \(H_\lambda\) and vanishing on the other eigenspaces. Therefore \(S_\lambda\) must commute with every operator in $\mathcal{B}(H_\lambda)$. Note that $\mathcal{B}(H_\lambda)'=\mathbb C I_{H_\lambda},$ then $S_\lambda=c_\lambda I_{H_\lambda}$ for some scalar \(c_\lambda\in\mathbb C\) and $\|S\|=\sup_{\lambda\in\Lambda}|c_\lambda|.$ Thus
\begin{equation*}
\{D\}''=\left\{\bigoplus_{\lambda\in\Lambda}c_\lambda I_{H_\lambda}:(c_\lambda)_{\lambda\in\Lambda}\in\ell^\infty(\Lambda)\right\}.
\end{equation*}
In particular, the map
\begin{equation*}
\Phi:\ell^\infty(\Lambda)\longrightarrow\{D\}'',\quad
\Phi(c)=\bigoplus_{\lambda\in\Lambda}c_\lambda I_{H_\lambda}
\end{equation*}
is an isometric algebra isomorphism. Furthermore, the same argument as in the proof of Theorem \ref{thm-diag} implies that  every weak $\ast$ continuous linear functional on \(\ell^\infty(\Lambda)\) is realized as a WOT matrix coefficient of \(\Phi(c)\), and then the WOT topology on \(\{D\}''\) is exactly $\sigma(\ell^\infty(\Lambda),\ell^1(\Lambda)).$

To be proceed, since for every polynomial $p\in\PP$,
\begin{equation*}
p(D)=\bigoplus_{\lambda\in\Lambda}p(\lambda)I_{H_\lambda}.
\end{equation*}
Therefore, denoting $\ell^{\infty}_D(\Lambda)=\left\{(p(\lambda))_{\lambda\in\Lambda}:p\in\PP\right\}$, we have
\begin{equation*}
\Phi^{-1}(\operatorname{alg}(D))=\ell^{\infty}_D(\Lambda).
\end{equation*}
Combining the topological relation as above, we obtain $\overline{\operatorname{alg}(D)}^{\,\mathrm{WOT}}=\{D\}''$
if and only if
\begin{equation}\label{equi-top-relation-DCP}
\overline{\ell^{\infty}_D(\Lambda)}^{\,\sigma(\ell^\infty(\Lambda),\ell^1(\Lambda))}
=\ell^\infty(\Lambda).
\end{equation}
Again by Lemma \ref{anni-dual}, \eqref{equi-top-relation-DCP} is equivalent to
${}_\perp\ell^{\infty}_D(\Lambda)=\{0\},$
where
\begin{equation*}
{}_\perp\ell^{\infty}_D(\Lambda)
=\left\{c\in\ell^1(\Lambda):\sum_{\lambda\in\Lambda}c_\lambda p(\lambda)=0\text{ for every }p\in\PP\right\}.
\end{equation*}
Since the monomials $(z^k)$ form a linear basis of the polynomial space, the above condition is equivalent to
\begin{equation*}
  \left\{c\in\ell^1(\Lambda):\sum_{\lambda\in\Lambda}c_\lambda \lambda^k=0\text{ for every }k\geq0\right\}=\{0\},
\end{equation*}
as we wanted.
\end{proof}

\medskip
\subsection{Examples}~
\medskip

To be anticipated, it is interesting for us to find that any real and bounded diagonal operator has the minimal commutant property on a separable Hilbert space if and only if its diagonal entries are pairwise distinct, as follows.

\begin{corollary}\label{cor-diag-real}
The real diagonal operator $D=\diag(d_n)$ has the minimal commutant property on a separable Hilbert space $H$ if and only if $(d_n)$ is bounded and pairwise distinct.
\end{corollary}

\begin{proof}
It is only needed to show the sufficiency by Theorem \ref{thm-diag}. For this, we argue by contradiction, and assume that there is some $c\in\ell^1\backslash\{0\}$ such that
\begin{equation*}
   \sum_{n=1}^{\infty}c_nd_n^k=0\ \mbox{ for all } k\geq0.
\end{equation*}
Set $R>0$ such that $d_n\in K=[-R,R]$ for all $n$. Consider the measure
\begin{equation*}
\mu=\sum_{n=1}^{\infty}c_n\delta_{d_n},
\end{equation*}
where $\delta_{d_n}$ is the Dirac measure on the single set $\{d_n\}$. Then $\mu$ is finite, because $\sum_{n=1}^{\infty}|c_n|<\infty$. Thus
\begin{equation*}
\int_K p(x)\,d\mu(x)=0\ \mbox{ for all } p\in \PP.
\end{equation*}
According to the Weierstrass Theorem, all polynomials is dense in the space of all continuous functions $C(K)$. So for any $f\in C(K)$, there are polynomials $p_j$ such that $\|p_j-f\|_{\infty}\to0$, and then
\begin{equation*}
\left|\int_K(f-p_j)\,d\mu\right|
\leq \|f-p_j\|_{\infty}\,|\mu|(K)\longrightarrow0.
\end{equation*}
Thus
\begin{equation*}
\int_K f\,d\mu=0,\qquad f\in C(K),
\end{equation*}
which in turn implies that $\mu=0$. Note that $c_n=\mu(\{d_n\})=0,$ $n\geq1,$ which leads to a contradiction.
\end{proof}

However, for the complex case, the following counterexample shows that there is a sequence of complex numbers inducing a bounded diagonal operator which fails to have a minimal commutant.

\begin{example}\label{ex-diag-ce}
Set
\begin{equation*}
b=\frac{1}{10},\quad r_m=1-\frac{1}{\sqrt{m}},\quad \omega_m=e^{2\pi i/m},\quad m\ge 2.
\end{equation*}
Note that $0<b<r_2<r_3<\cdots<1$. Let $d_1=0, c_1=1, d_2=b, c_2=-1.$
Define $N_m=2+\frac{(m-2)(m-1)}{2}\ (m\ge 2).$ For each integer $m\ge 2$, set
\begin{equation*}
d_{N_m+j+1}=r_m\,\omega_m^j,\quad
c_{N_m+j+1}=\frac{\lambda_m}{m}\,\omega_m^j,\quad j=0,1,\dots,m-1,
\end{equation*}
where $\lambda_m\in\mathbb C$ is defined recursively by
\begin{equation}\label{lam-rec-def}
\lambda_m =\left(\frac{b}{r_m}\right)^{m-1}-\sum_{\substack{\ell\mid m\\ 2\le \ell<m}}\lambda_\ell\left(\frac{r_\ell}{r_m}\right)^{m-1},
\end{equation}
where the sum runs over proper divisors $\ell$ of $m$, and the empty sum is zero. 

Then $(d_n)$ is bounded and pairwise distinct, and $(c_n)\in \ell^1\backslash\{0\}$ such that
$\sum_{n=1}^{\infty} c_n d_n^k=0$ for all integers $k\ge 0.$
\end{example}

\begin{proof}
We verify those properties in turn.

(1) The sequence $(d_n)$ is pairwise distinct and satisfies $|d_n|<1$ for all $n$. Indeed, for each group index $m\ge 2$, $|d_{N_m+j+1}|=r_m=1-1/\sqrt{m}<1$. For a fixed group $m$, $\{\omega_m^j\}_{j=0}^{m-1}$ are distinct, so the $m$ points in one group are pairwise distinct. The sequence $(r_m)_{m\ge 2}$ is strictly increasing, hence points from different groups have different moduli; these moduli are also strictly larger than $b$, so they cannot coincide with $d_1$ or $d_2$. Thus all $d_n$ are pairwise distinct and bounded by 1.

(2) For every $m\ge 2$ and every integer $k\ge 0$,
\begin{equation*}
\sum_{j=0}^{m-1}\frac{\lambda_m}{m}\,\omega_m^j\,(r_m\omega_m^j)^k =
\begin{cases}
 \lambda_m\, r_m^{k}, & m\mid (k+1),\\
 0, & m\nmid (k+1).
 \end{cases}
\end{equation*}
Indeed, rearrange the sum
\begin{equation*}
\sum_{j=0}^{m-1}\frac{\lambda_m}{m}\omega_m^j (r_m\omega_m^j)^k
=\frac{\lambda_m r_m^k}{m}\sum_{j=0}^{m-1}\omega_m^{(k+1)j}.
\end{equation*}
The geometric-series identity for roots of unity gives
\begin{equation*}
\sum_{j=0}^{m-1}\omega_m^{(k+1)j}=
\begin{cases}
m, & m\mid(k+1),\\
0, & m\nmid(k+1).
\end{cases}
\end{equation*}
So the desired result is obtained.

(3) $\displaystyle\sum_{n=1}^{\infty} c_n d_n^k = 0$ for all $k\ge 0$.
Indeed, for $k=0$, each group satisfies
\begin{equation*}
 \sum_{j=0}^{m-1} c_{N_m+j+1}= \frac{\lambda_m}{m}\sum_{j=0}^{m-1}\omega_m^j=0,
\end{equation*}
and the initial point contribution is $c_1+c_2=1-1=0$. Hence the total sum vanishes.

Now let $k\ge 1$ and set $M:=k+1\ge 2$. The total moment splits into the initial point contribution plus contributions from all groups,
\begin{equation*}
\sum_{n=1}^\infty c_n d_n^k=c_1 d_1^k + c_2 d_2^k+\sum_{m\ge 2}\sum_{j=0}^{m-1} c_{N_m+j+1}d_{N_m+j+1}^k.
\end{equation*}
Since $d_1=0$, $c_1 d_1^k=0$, and $c_2 d_2^k = -b^{k}=-b^{M-1}$. By (2), only groups with $m\mid M$ including non-zero terms, so
\begin{equation*}
\sum_{n=1}^\infty c_n d_n^k
=-b^{M-1}+\sum_{\substack{m\mid M\\ m\ge 2}}\lambda_m r_m^{M-1}.
\end{equation*}
From recursion \eqref{lam-rec-def}, substitute $m=M$ and multiply both sides by $r_M^{M-1}$:
\begin{equation*}
\lambda_M r_M^{M-1}=b^{M-1}-\sum_{\substack{\ell\mid M\\ 2\le \ell<M}}\lambda_\ell r_\ell^{M-1}.
\end{equation*}
That is to say,
$-b^{M-1}+\sum_{\substack{\ell\mid M\\ \ell\ge 2}}\lambda_\ell r_\ell^{M-1}=0.$
Thus the $k$-th moment equals zero. For fixed $k$, only finitely many divisors $m$ of $M=k+1$ appear, so the computation reduces to a finite sum.

(4) $\displaystyle\sum_{m=2}^{\infty}|\lambda_m|<\infty$.
For this, let $\rho = b/r_2\in(0,1)$ and $\gamma=\frac{\sqrt{2}-1}{2}>0$. From \eqref{lam-rec-def},
\begin{equation*}
|\lambda_m|\le \rho^{m-1}
+\sum_{\substack{\ell\mid m\\ 2\le \ell<m}}|\lambda_\ell|\left(\frac{r_\ell}{r_m}\right)^{m-1}.
\end{equation*}
If $\ell$ is a proper divisor of $m$, then $\ell\le m/2$. Using $r_k=1-k^{-1/2}$, one obtains the estimate
\begin{equation*}
\left(\frac{r_\ell}{r_m}\right)^{m-1}\le \exp\big(-\gamma\sqrt{m}\big).
\end{equation*}
Define $S_m=\sum_{\ell=2}^{m}|\lambda_\ell|$. Then
$S_m \le \big(1+e^{-\gamma\sqrt{m}}\big)S_{m-1}+\rho^{m-1}.$
Iteration gives
\begin{equation*}
S_m \le \sum_{j=2}^{m}\rho^{j-1}\prod_{k=j+1}^{m}\big(1+e^{-\gamma\sqrt{k}}\big).
\end{equation*}
Since $\sum_{t=2}^{\infty}e^{-\gamma\sqrt{t}}<\infty$, the infinite product
$P=\prod_{k=2}^{\infty}\big(1+e^{-\gamma\sqrt{k}}\big)$
converges to a finite constant. Hence
$S_m \le P\sum_{j=2}^{\infty}\rho^{j-1}= P\,\frac{\rho}{1-\rho}<\infty.$ So $\sum_{m=2}^{\infty}|\lambda_m|=\lim_{m\to\infty}S_m<\infty$. Note that $c_1=1$ and
\begin{equation*}
\sum_{n=1}^{\infty}|c_n|
=|c_1|+|c_2|+\sum_{m=2}^{\infty}\sum_{j=0}^{m-1}\left|\frac{\lambda_m}{m}\omega_m^j\right|
=2+\sum_{m=2}^{\infty}|\lambda_m|<\infty.
\end{equation*}
So $(c_n)\in \ell^1\backslash\{0\}$.
\end{proof}

As is illustrated in Example \ref{ex-diag-ce}, it is not easy to verify condition (b) of Theorem \ref{thm-diag}. The following highly non-trivial examples illustrate again this phenomenon. For this, we recall some results about atomic measures and Poisson kernels.

For $0<r<1$, let $P_r$ denote the Poisson kernel
$P_r(e^{it})=\frac{1-r^2}{|e^{it}-r|^2}.$
For a finite measure $\nu$ on $\TT$, denote by
\begin{equation*}
(P_r*\nu)(z):=\int_{\mathbb T}P_r(z\overline w)\,d\nu(w), \quad z\in\mathbb T.
\end{equation*}
By the Fourier expansion of the Poisson kernel, we have
\begin{equation*}
P_r(z\overline w)=\sum_{n\in\mathbb Z}r^{|n|}z^n\overline w^{\,n},
\end{equation*}
where the series converges absolutely and uniformly. Then we may interchange the order of summation and integration, and obtain
\begin{equation*}
(P_r*\nu)(z)=\sum_{n\in\mathbb Z}r^{|n|}
z^n\int_{\mathbb T}\overline w^{\,n}\,d\nu(w).
\end{equation*}
Moreover, since $P_r*\nu$ is a continuous function on $\mathbb T$, then the measure $d\sigma(z)=(P_r*\nu)(z)\,dm(z)$
is continuous absolutely with respect to normalized Lebesgue measure $m$. On the other hand, for any $k\geq0$, by Fubini's Theorem,
\begin{equation*}
\int_{\mathbb T}z^k(P_r*\nu)(z)\,dm(z)
=\int_{\mathbb T}
\left[\int_{\mathbb T}
z^kP_r(z\overline w)\,dm(z)\right]d\nu(w).
\end{equation*}
Note that
\begin{equation*}
\int_{\mathbb T}z^kP_r(z\overline w)\,dm(z)
=\sum_{n\in\mathbb Z}r^{|n|}\overline w^{\,n}\int_{\mathbb T}z^{k+n}\,dm(z)=r^k w^k,
\end{equation*}
which in turn implies that
\begin{equation*}
\int_{\mathbb T}z^k(P_r*\nu)(z)\,dm(z)
=r^k\int_{\mathbb T}w^k\,d\nu(w), \quad k\geq0.
\end{equation*}
In other words,
\begin{equation*}
\int_{\mathbb T}z^k P_r*\nu(z)\,dm(z)=r^k\int_{\mathbb T}z^k\,d\nu(z),
\quad k\geq0.
\end{equation*}
Also recall that F. and M. Riesz Theorem: for a finite measure $\mu$,
\begin{equation*}
\int_{\mathbb T}z^k\,d\mu(z)=0 \mbox{ for all } k\geq0 \quad \Longrightarrow \quad \mu\ll m,
\end{equation*}
(for example see \cite[Theorem 17.13]{Rudin87}).

\begin{corollary}\label{cor:J-one}
Let $\Omega$ be the bounded area enclosed by some Jordan curve $\Gamma$. Let the sequence $(d_n)$ be pairwise distinct. If $(d_n)=(x_p)\cup(y_q)$ with $(x_p)$ dense in $\Gamma$, and $(y_q)$ converging to some $a\in\Omega$, then the corresponding diagonal operator $D=\diag\{d_n\}$ has the minimal commutant property on a separable Hilbert space.
\end{corollary}

\begin{proof}
Given any $c\in\ell^{1}$ such that
\begin{equation*}
  \sum_{n=1}^{\infty}c_n d_n^k=0\quad \mbox{for all } k\geq0,
\end{equation*}
we consider the following two finite measures
\begin{equation*}
\nu=\sum_p\alpha_p\delta_{x_p},\quad
\rho=\sum_q\beta_q\delta_{y_q},
\end{equation*}
where $c=(\alpha_p)\cup(\beta_q)$ with respect to the corresponding decomposition of $(d_n)=(x_p)\cup(y_q)$. Then
\begin{equation*}
\int_\Gamma z^m\,d\nu(z)+\int_\Omega z^m\,d\rho(z)=0,\quad m\geq0.
\end{equation*}
Note that $\rho$ is purely atomic and its support
$\supp(\rho)\subseteq\{y_q:q\geq1\}\cup\{a\}$,
contained in the interior $\Omega^\circ$. And
\begin{equation*}
  \dist(\supp(\rho),\Gamma)>0.
\end{equation*}

For $y_q\in\Omega$, denote by $\omega_{y_q}$ the harmonic measure with respect to $y_q$ on the domain $\Omega$. Then for any function $f$ analytic inside $\Omega$ and continuous on  $\overline\Omega$, we have
\begin{equation*}
f(y_q)=\int_\Gamma f(z)\,d\omega_{y_q}(z).
\end{equation*}
Now define $\sigma=\sum_q\beta_q\omega_{y_q}.$ Clearly, $\sigma$ is a finite measure on $\Gamma$ by the fact that $\sum_1|\beta_q|\leq \|c\|_1<\infty$, and
\begin{equation}\label{mea-inter-bound}
\int_\Omega f\,d\rho=\int_\Gamma f\,d\sigma.
\end{equation}
So
\begin{equation}\label{mea-inter-bound=0}
\int_\Gamma z^m\,d(\nu+\sigma)(z)=0,\quad m\geq0
\end{equation}

On the other hand, by the Riemann Mapping Theorem and Carath\'{e}odory Theorem, there is a conformal map
\begin{equation*}
\varphi:\Omega\longrightarrow\mathbb D,
\end{equation*}
which can be continuous extension to its boundary. Denote by $m$ the normalized Lebesgue measure on $\TT$. For any Borel subset $E\subseteq\Gamma$,
\begin{equation*}
\omega_z(E)=\omega_u\bigl(\varphi^{-1}(E)\bigr)\ \ \mbox{with }  u=\varphi^{-1}(z).
\end{equation*}
Note that
\begin{equation*}
\omega_u(F)=\int_F\frac{1-|u|^2}{|\zeta-u|^2}\,dm(\zeta),\quad F\subseteq\TT,
\end{equation*}
which in turn implies that
\begin{equation*}
\omega_z(E)=\int_{\varphi^{-1}(E)}\frac{1-|u|^2}{|\zeta-u|^2}\,dm(\zeta).
\end{equation*}
So
\begin{equation*}
\omega_a(E)=\int_{\varphi^{-1}(E)}\frac{1-|v|^2}{|\zeta-v|^2}\,dm(\zeta).
\end{equation*}
Moreover, it is easy to see that
\begin{equation*}
\frac{1-|u|}{1+|u|}\leq P_u(\zeta)\leq\frac{1+|u|}{1-|u|}.
\end{equation*}
Then
\begin{equation*}
C^{-1}\omega_a(E)\leq\omega_z(E)\leq C\omega_a(E),
\end{equation*}
for some $C>0$. (For more details, see \cite[Chapter I]{GM}.) Then all the harmonic measure with respect to the interior point $z_j\in\Omega$ is absolutely continuous to the harmonic measure $\omega_a$, and thus $\sigma\ll\omega_a,$ which in turn implies that
\begin{equation*}
\nu+\sigma\ll\omega_a
\end{equation*}
by the F. and M. Riesz Theorem. Therefore, we have $\nu=(\nu+\sigma)-\sigma\ll\omega_a.$ Since  $\nu$ is purely atomic but $\omega_a$ has no atomic, then $\nu=0$ and $x_p=0$ for $p\geq1$. That, together with \eqref{mea-inter-bound} and \eqref{mea-inter-bound=0}, implies
\begin{equation*}
  \int_\Omega z^m\,d\rho(z)=0, \quad m\geq0.
\end{equation*}
Thus for any polynomial $P$, we have
\begin{equation}\label{poly=0}
  \int_\Omega P(z)\,d\rho(z)=0.
\end{equation}

Consider the compact set
\begin{equation*}
K_0=\{y_q:q\geq1\}\cup\{a\}\subseteq\Omega,
\end{equation*}
with its interior $K_0^\circ=\varnothing.$ So
$\mathbb C\setminus K_0$ is connected, and then by Mergelyan's Theorem (for example see \cite{Ga}),
\begin{equation*}
\overline{\mathcal P(K_0)}=C(K_0).
\end{equation*}
Together with \eqref{poly=0}, it follows that
\begin{equation*}
  \int_{K_0} f(z)\,d\rho(z)=0
\end{equation*}
for any $f\in C(K_0)$. So $\rho=0$ by the Riesz Representation Theorem, and then $x_q=0$. Now the proof is complete.
\end{proof}

\medskip
\subsection{Sufficient criteria}~
\medskip

Inspired by \cite[Theorem 5.3]{BY}, we provide some easily verifiable criteria via the geometric properties of the spectrum of the diagonal operator, by using the complex-polynomial approximation. 
As usual, denote by $A(\Omega)$ the algebra consists of all functions analytic inside a bounded domain $\Omega\subseteq\CC$ and continuous extension to its boundary $\partial\Omega$, with the uniform closure of the polynomials. We shall use the following classical extension theorem (for example, see \cite{Carleson} and \cite{Rudin56}).

\begin{theorem}[Rudin-Carleson]
Let $E\subseteq\TT$ be a closed set of Lebesgue measure zero and let $f:E\rightarrow\CC$ be a continuous function. Then there exists $g\in A(\DD)$ such that $g|_E=f$ and $\|g\|_\DD\leq \|g\|_E$.
\end{theorem}

As an application, the following, \cite[Theorem 5.3]{BY}, holds.

\begin{proposition}\label{thm:diagonal}
Let $H$ be a Hilbert space with orthonormal basis $(e_n)_{n\in\mathbb N}$, let $(d_n)$
be a bounded sequence of complex numbers and let $D$ be the diagonal operator defined by $D(e_n)=d_n e_n$. Assume that, for all bounded sequences of complex numbers $(\mu_n)$, there exists $M>0$ such that, for all $\varepsilon>0,$ for all $N\geq 1,$ there exists $P\in\PP$ such that
 $$\left\{
 \begin{array}{rcl}
  |P(d_n)-\mu_n|&\leq&\varepsilon,\ n=1,\dots,N\\
  \sup_{n\in\NN}|P(d_n)|&\leq&M.
 \end{array}\right.$$
Then $D$ has a minimal commutant.
\end{proposition}

However, the converse of Proposition \ref{thm:diagonal} does not hold, in general. Indeed, as in Corollary \ref{cor:J-one}, let $\Omega$ be the bounded area enclosed by some Jordan curve $\Gamma$, and let the pairwise distinct  sequence $(d_n)$ satisfy that $(d_n)=(d_{2k-1})\cup(d_{2k})$ with $\Lambda_1=(d_{2k-1})$ dense in $\Gamma$ and $\Lambda_2=(d_{2k})$ converging to some $a\in\Omega$. By Corollary \ref{cor:J-one}, $D=\diag(d_n)$ has a minimal commutant. But the condition of Proposition \ref{thm:diagonal} dose not hold. To end this, consider the following diagonal operator $Te_{2k-1}=0$; and meanwhile $Te_{2k}=(-1)^{k}e_{2k}$. Clearly, $T\in\{D\}'$.  If this is not the case, then for the corresponding sequence $(\mu_n)$ with $\mu_{2k-1}=0$ and $\mu_{2k}=(-1)^{k}$, there exists $M>0$ such that, for all $\varepsilon>0,$
 for all $N\geq 1,$ there exists $P\in\PP$ such that
 $$\left\{
 \begin{array}{rcl}
  |P(d_n)-\mu_n|&\leq&\varepsilon,\ n=1,\dots,N\\
  \sup_{n\in\NN}|P(d_n)|&\leq&M.
 \end{array}\right.$$
Then one can choose a sequence of polynomials $(P_j)\subseteq \PP$ satisfying the above conditions, and thus $P_j(D)$ converges weakly to $T$ by Proposition \ref{thm:diagonal}. By the principle of uniform boundedness, there is some constant $C>0$ such that
\begin{equation*}
  \sup_{j}\|P_j(D)\|\leq C.
\end{equation*}
By the density of $\Lambda_1$ in $\Gamma$, we have
\begin{equation*}
  \sup_{z\in\Gamma}|P_j(z)|=\sup_{d_{2k-1}\in\Lambda_1}|P_j(d_{2k-1})|\leq C,
\end{equation*}
which in turn implies that
\begin{equation*}
  \sup_{z\in\overline{\Omega}}|P_j(z)|\leq C
\end{equation*}
due to the maximum modulus principle. On the other hand, for the subsequence $\Lambda_2$ converging to $a\in\Omega$, we have
\begin{equation*}
  P_j(d_{2k})\rightarrow(-1)^{k}.
\end{equation*}
By Montel's theorem, there is a subsequence of $(P_j)$, also denoted by $(P_j)$, such that
$P_j$ converges uniformly to an analytic function $f$ on any compact subset of $\Omega$. Particularly,
\begin{equation*}
  P_j(y_{q})\rightarrow f(y_{q}).
\end{equation*}
So $f(d_{2k})=(-1)^{k}$, which leads to a contradiction by the hypotheses $\Lambda_2$ converging to $a\in\Omega$.

On the other hand, inspired by Proposition \ref{thm:diagonal}, the following sufficient criteria holds, generalizing \cite[Proposition 5.4 (1)]{BY}.

\begin{theorem}\label{thm:sepc-mc}
Let $(d_n)$ be bounded and pairwise distinct.
Let $H$ be a Hilbert space with orthonormal basis $(e_n)_{n\in\mathbb N}$,  and let $D$ be the diagonal operator defined by $D(e_n)=d_n e_n$.
\begin{itemize}
  \item[(1)] Let the complement $\CC\backslash \overline{\{d_n:n\geq1\}}$  be connected. If $(d_n)$ can be decomposed into the union of finitely many pairwise disjoint convergent subsequences, then $D=\diag(d_n)$ has a minimal commutant.
  \item[(2)] If $(d_n)$ is contained in some Jordan curve, then $D$ has a minimal commutant.
  \item[(3)] Let $(d_n)=S_1\cup\cdots\cup S_m$ for some pairwise disjoint sets $S_1,\ldots,S_m$, and let $K_j=\overline{S_j}$ be pairwise disjoint compact totally disconnected sets. Denote $K:=\bigcup_{j=1}^m K_j$ and let $\CC\setminus K$ be connected. Then $D$ has a minimal commutant.
\end{itemize}
\end{theorem}

\begin{proof}
Since $(d_n)$ is pairwise distinct, then $T\in\{D\}'$ if and only if $T$ is of the form $T=\diag(\mu_n)$ for some bounded sequence $(\mu_n)$. Set $M=\sup_{n\ge1}|\mu_n|=\|T\|<+\infty.$

(1) can be obtained immediately from \cite[Proposition 5.4 (1)]{BY} and Corollary \ref{sum-normal}.

(2) Let $(d_n)\subseteq\Gamma$, a Jordan curve,  and let $\Omega$ be the bounded area enclosed by the curve $\Gamma$. By the Riemann mapping theorem, there exists a biholomorphic mapping $\phi:\DD\rightarrow\Omega$. Since $\Gamma$ is a Jordan curve, by Carath\'{e}odory boundary extension theorem, $\phi$ has a continuous extension to the closed unit disk, and $\phi(\TT)=\Gamma$. Let $\lambda_n=\phi^{-1}(d_n)$. Then $(\lambda_n)$ is pairwise distinct.

Let $T\in\{D\}'$. There exists a sequence $(\mu_n)\in\ell^{\infty}$, bounded by $\|T\|$, such that $Te_n=\mu_n e_n$ for all $n\in\NN$. Let $N\in \mathbb{N}$ be fixed. Consider the finite set
$E_N=\{d_1,d_2, \dots, d_N\}$, and the continuous function
$g_N \colon E_N \to \mathbb{C}$, defined by the expression $g_N(d_k)=\mu_k$, for every $k=1,2,\dots,N$. It follows from Rudin-Carleson's Theorem that there exists a function $g_N\in A(\mathbb{D})$
such that $g_N(d_k)=\mu_k$ for every $k=1,2,\dots,N$, and such that
$\|g_N\|_{\mathbb{D}} \le \|(d_k)\|_{\infty}$.

Next, consider $h_N(z)=g_N(\phi^{-1}(z))$, and then $h_N\in A(\Omega)$ the algebra with the uniform closure of the polynomials on $\Omega$, because $\phi^{-1}$ has a continuous extension to $\Omega$. Clearly, $h_N(d_n)=g_N(\phi^{-1}(d_n))=\mu_k$. Meanwhile, we have
\begin{equation*}
  \|h_N\|_{\overline{\Omega}}=\|g_N\|_{\overline{\DD}}\leq\|(\mu_k)\|_{\infty}.
\end{equation*}

Note that the complement $\overline{\Omega}^{c}$ is connected. Then by Mergelyan’s Theorem, for $N\geq 1,$ there exists $P_N\in\PP$ be such that
\begin{equation*}
  \|P_N-h_N\|_{\overline{\Omega}}<\frac{1}{N}.
\end{equation*}
Thus, $\sup_{n\in\NN}|P_N(d_n)|\leq M$ and $|P_N(d_n)-\mu_n|\leq 1/N$ for $n=1,\dots,N$ and $M>0$. So the desired result can be immediately obtained by Proposition \ref{thm:diagonal}.

(3) Fix a positive integer $N$. For each $j=1,\ldots,m$, define
$F_{j,N}=\{d_k:1\le k\le N,\ \lambda_k\in S_j\}.$ Since each $K_j$ is compact and totally disconnected, it is a zero dimensional compact Hausdorff space. Hence, for the finite set $F_{j,N}$, we can find pairwise disjoint relatively open and closed subsets
$V_{j,1},\ldots,V_{j,r_j}\subseteq K_j$ such that $F_{j,N}=\{d_{k_1},\ldots,d_{k_{r_j}}\},$ $\lambda_{k_l}\in V_{j,l}.$ Define
$$f_{j,N}(z)=\begin{cases}\mu_{k_l},&z\in V_{j,l},\quad 1\le l\le r_j,\\
0,&z\in K_j\setminus\displaystyle\cup_{l=1}^{r_j}V_{j,l}.\end{cases}$$
Clearly, $f_{j,N}(d_k)=\mu_k,$ $d_k\in F_{j,N}$. Because the sets $V_{j,l}$ are clopen, then $f_{j,N}$ is continuous on $K_j$ with $\|f_{j,N}\|_{K_j}\le M.$

Since each $K_j$ is compact and totally disconnected,
$\mathbb C\setminus K_j$
is connected. By Mergelyan’s theorem, there exists a polynomial $q_{j,N}\in\PP$ satisfying
\begin{equation*}
\sup_{z\in K_j}|q_{j,N}(z)-f_{j,N}(z)|<\frac{1}{N},
\end{equation*}
which in turn implies that $\sup_{z\in K_j}|q_{j,N}(z)|\le M+\frac{1}{N}.$

Since $K_1,\ldots,K_m$ are pairwise disjoint compact sets, we may choose pairwise disjoint open neighbourhoods
$U_1,\ldots,U_m$ such that $K_j\subseteq U_j.$
Set $U=\cup_{j=1}^m U_j.$ Define, on $U$,
$h_N(z)=q_{j,N}(z),$ $z\in U_j.$ Then $h_N$ is analytic on $U$, due to the fact $U_j$ are pairwise disjoint. Using the hypothesis
$\CC\setminus K$ is connected, it follows from Runge’s theorem that there is a polynomial $P_N$ such that
$$\sup_{z\in K}|P_N(z)-h_N(z)|<\frac{1}{N}$$
and meanwhile
$$\sup_{z\in K}|P_N(z)|\le M+\frac{2}{N}.$$
Thus $\sup_{N\ge1}\|P_N(D)\|\le M+1.$
On the other hand, for $1\le k\le N$, suppose $d_k\in S_j$. Then
\begin{equation*}
|P_N(d_k)-\mu_k|\le |P_N(d_k)-q_{j,N}(d_k)|+|q_{j,N}(d_k)-f_{j,N}(d_k)|<\frac{2}{N},
\end{equation*}
i.e.,
\begin{equation*}
|P_N(d_k)-\mu_k|<\frac{2}{N},\ 1\le k\le N.
\end{equation*}

Finally, for arbitrary $x=\sum_{n=1}^{\infty}x_ne_n\in H$, we have
$$\|P_N(D)x-Tx\|^2=\sum_{n=1}^{\infty}|P_N(d_n)-\mu_n|^2|x_n|^2\longrightarrow0.$$
Therefore $P_N(D)\rightarrow T$ in the operator norm topology, which implies that $D$ has a minimal commutant.
\end{proof}

\section{Tiny commutant of diagonal operators}\label{sec:tiny-com-diag}

This section gives similar criteria when a diagonal operator has the tiny minimal or double commutant property on a separable Hilbert space. For this, given a compact set $K\subseteq\CC$, denote
\begin{equation*}
\Rat(K)=\left\{f: \mbox{ rational functions } f \mbox{ with poles outside } K\right\}.
\end{equation*}

\begin{theorem}\label{thm-diag-tiny}
Let $D=\diag(d_n)$ be a bounded diagonal operator on a separable Hilbert space $H$, and denote $K=\overline{\{d_n:n\geq1\}}$. Then the following are equivalent:
\begin{itemize}
  \item[(i)] $D$ has the tiny commutant property.
  \item[(ii)] $(d_n)$ is pairwise distinct, and \begin{equation}\label{eq:tiny-com-ann1}
  \left\{c=(c_n)\in\ell^1: \sum_{n=1}^{\infty}\frac{c_n}{z-d_n}=0\ \mbox{ for all } z\in\CC\setminus K\right\}=\{0\}.
\end{equation}
  \item[(iii)] $(d_n)$ is pairwise distinct, and there is no nonzero finite purely atomic complex measure
   \begin{equation*}
   \mu=\sum_{n=1}^{\infty}c_n\delta_{d_n}\ \mbox{ with } c=(c_n)\in\ell^1
   \end{equation*}
   whose Cauchy transform vanishes identically on \(\mathbb C\setminus K\).
\end{itemize}
\end{theorem}

\begin{proof}
(i)$\Leftrightarrow$(ii) If $D$ has a tiny commutant, then $(d_n)$ is pairwise distinct. Indeed, if this is not the case, then there are some $j\neq m$ such that $d_j=d_m$. Consider the operator $T=e_j\otimes e_m$. It is easy to see that $Te_m=e_j$ and $Te_k=0$ for all $k\neq m$, and then $T\in\{D\}'$, $\langle Te_m,e_j\rangle=1$. On the other hand, for any $f\in\Rat(K)$, $\langle f(D)e_m,e_j\rangle=0$. So $X\in\{D\}'\backslash\overline{\alg^{a}\{D\}}^{\,WOT}$, which leads to a contradiction. In the sequel, we can always assume that $(d_n)$ is pairwise distinct. Thus
\begin{equation*}
\{D\}'=\left\{\diag(\mu_n):(\mu_n)\in\ell^\infty\right\}\subseteq \mathcal{B}(H).
\end{equation*}
Now similar to the proof of Theorem \ref{thm-diag}, $D$ has a tiny commutant if and only if
\begin{equation}\label{eq:tiny-com-ann0}
\left\{c=(c_n)\in\ell^1:\sum_{n=1}^{\infty}c_nf(d_n)=0,\ \mbox{ for all } f\in\Rat(K)\right\}=\{0\}.
\end{equation}

Furthermore, we claim that \eqref{eq:tiny-com-ann0} is equivalent to \eqref{eq:tiny-com-ann1}. In fact, given any $c=(c_n)\in\ell^1$, denote $F_c(z)=\sum_{n=1}^{\infty}\frac{c_n}{z-d_n}$ for short, which is analytic in $\CC\setminus K$. If \eqref{eq:tiny-com-ann0} holds, taking $f(w)=\frac1{z-w},$ $z\in\CC\setminus K,$ we have $F_c(z)=0$ $(z\in\CC\setminus K)$. Conversely, assume that $F_c(z)=0$ $(z\in\CC\setminus K)$. For any $|z|>\sup_n|d_n|$, consider its Laurent series expansion
\begin{equation*}
F_c(z)=\sum_{n=1}^{\infty}c_n\sum_{k=0}^{\infty}\frac{d_n^k}{z^{k+1}}=\sum_{k=0}^{\infty}\frac{\sum_{n=1}^{\infty}c_nd_n^k}{z^{k+1}}.
\end{equation*}
Then $\sum_{n=1}^{\infty}c_nd_n^k=0,$ $k=0,1,2,\ldots.$ Moreover, given any $z\in\CC\setminus K$, by considering the $m$-th derivative of $\frac1{z-w}$ at $w=d_n$, we have
\begin{equation*}
\sum_{n=1}^{\infty}\frac{c_n}{(z-d_n)^m}=0, \ \mbox{ for all } z\in\CC\setminus K.
\end{equation*}
Since every rational function can be written as the sum of a polynomial and finitely many principal parts, then
\begin{equation*}
\sum_{n=1}^{\infty}c_nf(d_n)=0,\ \mbox{ for all } f\in\Rat(K).
\end{equation*}
So \eqref{eq:tiny-com-ann0} is equivalent to \eqref{eq:tiny-com-ann1}.

(ii)$\Leftrightarrow$(iii) Assume that $(d_n)$ is pairwise distinct. Given any $c=(c_n)\in\ell^1$, consider the measure
\begin{equation*}
\mu=\sum_{n=1}^{\infty}c_n\delta_{d_n},
\end{equation*}
and then $\|\mu\|=\sum_{n=1}^{\infty}|c_n|<\infty.$ Clearly, $\mu$ is a finite purely atomic complex measure, and its Cauchy transform is
\begin{equation*}
\mathcal C\mu(z)=\int_K\frac{1}{z-w}\,d\mu(w)=\sum_{n=1}^{\infty}\frac{c_n}{z-d_n},\ \mbox{ for all } z\in\CC\setminus K.
\end{equation*}
So the desired result is obtained.
\end{proof}

Now we can also provide some easily verifiable criteria via the geometric properties of the spectrum of the diagonal operator.

\begin{corollary}
Given a sequence $(d_n)$ consist of pairwise distinct complex numbers, denote $K=\overline{\{d_n:n\geq1\}}$, and let $D=\diag(d_n)$ be a bounded diagonal operator on a separable Hilbert space $H$.
\begin{itemize}
  \item[(i)] If $K$ has planar Lebesgue measure zero, then $D$ has a tiny commutant.
  \item[(ii)] If $\CC\setminus K$ is connected, and the interior $K^{\circ}=\varnothing$, then $D$ has a tiny commutant.
\end{itemize}
\end{corollary}

\begin{proof}
(i) We denote by $m_2$ the planar Lebesgue measure and assume $m_2(K)=0$. Given any $c=(c_n)\in\ell^1$, consider the finite measure
\begin{equation*}
\mu=\sum_{n=1}^{\infty}c_n\delta_n.
\end{equation*}
If $\mathcal C\mu\equiv0$ for all $z\in\CC\setminus K$, then $\mathcal C\mu=0$ a.e. on $\CC$, and $\mathcal C\mu\in L^1_{\mathrm{loc}}(\CC)$, the space of locally integrable functions. So the corresponding distribution of $\mathcal C\mu$ is the zero distribution. Moreover, recall the Cauchy kernel distribution identity
\begin{equation*}
\bar\partial\left(\frac1{\pi z}\right)=\delta_0,
\end{equation*}
where $\delta_0$ is the Dirac measure at the point $0$. Then $\bar\partial\mathcal C\mu=\pi\mu,$ and thus $\mu=0$. So the desired result holds immediately from Theorem \ref{thm-diag-tiny}.

(ii) Similarly, given any $c=(c_n)\in\ell^1$, consider the corresponding finite measure $\mu$ and assume $\mathcal C\mu\equiv0\text{ on }\CC\setminus K$. Then
\begin{equation*}
\mathcal C\mu(z)=\int_K\frac{d\mu(w)}{z-w}=0,\quad z\in\CC\setminus K.
\end{equation*}
Note that there is some $R>0$ such that $K\subseteq\{z:|z|<R\}$ because $K$ is compact. Moreover, for any $|z|>R$,
\begin{equation*}
\frac1{z-w}=\frac1z\frac1{1-w/z}=\sum_{n=0}^{\infty}\frac{w^n}{z^{n+1}},
\end{equation*}
and then
\begin{equation*}
\mathcal C\mu(z)=\sum_{n=0}^{\infty}\frac{1}{z^{n+1}}\int_Kw^n\,d\mu(w).
\end{equation*}
Using the uniqueness of the Laurent series expansion,
\begin{equation*}
\int_K w^n\,d\mu(w)=0,\quad n=0,1,2,\ldots,
\end{equation*}
which in turn implies that $\int_Kp(w)\,d\mu(w)=0$ for any $p\in\PP$. On the other hand, since $\CC\setminus K$ is connected, and $ K^\circ=\varnothing$, then by the Mergelyan's Theorem, for any $f\in C(K)$, there is a sequence $(p_j)\subseteq\PP$ such that $\|p_j-f\|_K\rightarrow0.$ So
\begin{equation*}
\left|\int_K f\,d\mu\right|=\left|\int_K(f-p_j)\,d\mu\right|\leq \|f-p_j\|_K\,|\mu|(K)\rightarrow0,
\end{equation*}
in other words,
\begin{equation*}
\int_Kf\,d\mu=0,\quad f\in C(K).
\end{equation*}
Therefore, $\mu=0$ by the Riesz Representation Theorem. Again by Theorem \ref{thm-diag-tiny}, $D$ has a tiny commutant.
\end{proof}

There exists a diagonal operator who has a tiny commutant but fails to have the minimal commutant property as follows.

\begin{example}
Let $(d_n)$ come from Example \ref{ex-diag-ce}. Then $D=\diag(d_n)$ has a tiny commutant on a separable Hilbert space $H$.
\end{example}

\begin{proof}
We use the same notation as in Example \ref{ex-diag-ce}. Note that
\begin{equation*}
\overline{\bigcup_{m\ge2}\{r_m\omega_m^j:0\le j<m\}}=\bigcup_{m\ge2}\{r_m\omega_m^j:0\le j<m\}\cup\TT.
\end{equation*}
Then
\begin{equation*}
\sigma(D)=\{0,b\}\cup\bigcup_{m\ge2}r_m\UU_m\cup\TT,
\end{equation*}
where $\UU_m=\{\omega_m^j:0\le j<m\}.$ So $m_2(\sigma(D))=0$, and then $D=\diag(d_n)$ has a tiny commutant by Theorem \ref{thm-diag-tiny}.
\end{proof}

At the end, the tiny double commutant property can be characterized by combining Theorem \ref{thm-diag-DCP} and Theorem \ref{thm-diag-tiny}.

\begin{theorem}\label{thm-diag-tiny-double}
Let $D=\diag(d_n)$ be a bounded diagonal operator on a separable Hilbert space $H$. Let $\Lambda=\{d_n:n\geq 1\}$ be the set of its distinct diagonal entries and $K=\overline{\Lambda}$. Then $D$ has the tiny double commutant property if and only if
\begin{equation*}
\left\{c=(c_\lambda)\in\ell^1(\Lambda):\sum_{\lambda\in\Lambda}
\frac{c_\lambda}{z-\lambda}=0,\text{ for every }z\in\CC\setminus K\right\}=\{0\}.
\end{equation*}
\end{theorem}

\begin{proof}
For each \(\lambda\in\Lambda\), define $I_\lambda=\{n\geq 1:d_n=\lambda\}$ and let $H_\lambda=\overline{\operatorname{span}}\{e_n:n\in I_\lambda\}.$
Then
\begin{equation*}
H=\bigoplus_{\lambda\in\Lambda}H_\lambda, \mbox{ and }
D=\bigoplus_{\lambda\in\Lambda}\lambda I_{H_\lambda}.
\end{equation*}
Thus, by the same argument as in the proof of Theorem \ref{thm-diag-DCP},
\begin{equation*}
\{D\}''=\left\{\bigoplus_{\lambda\in\Lambda}c_\lambda I_{H_\lambda}:(c_\lambda)_{\lambda\in\Lambda}\in\ell^\infty(\Lambda)\right\},
\end{equation*}
and then
\begin{equation*}
\left\{c=(c_\lambda)\in\ell^1(\Lambda):\sum_{\lambda\in\Lambda}c_\lambda f(\lambda)=0,\text{ for every } f\in\Rat(K)\right\}=\{0\}.
\end{equation*}
So the desired result can be obtained similar to the proof of Theorem \ref{thm-diag-tiny}.
\end{proof}

\section{Application to composition operators on Fock space}\label{sec:app-co}

As in \cite{BWY,BY}, some composition operators can be viewed as diagonal operators. Particularly, it follows from \cite[Theorem 3.1]{BY} that some composition operators $C_\varphi$ having the minimal commutant property restrict the canonical Jordan form of the matrix $\varphi'(0)$ at most one Jordan block of order $2$ on Hilbert spaces of holomorphic functions in several complex variables.

\medskip
\subsection{The case of diagonal matrixes}~

\medskip
To understand their commutants, we first consider the following highly non-trivial diagonal operators in Corollary \ref{ex-diag-mc}, by using a similar argument to Corollary \ref{cor:J-one}. For this, we need a generalization of the F. and M. Riesz Theorem.

\begin{lemma}\label{lem:ann}
Let $(r_j)$ be a strictly decreasing sequence of positive numbers, and let $(\mu_j)$ be finite purely atomic complex measures on $\TT$ satisfying $\sum_{j=0}^\infty\|\mu_j\|<\infty.$ If
\begin{equation}\label{lem:eq:ann}
\sum_{j=1}^\infty r_j^k\int_{\TT}z^k\,d\mu_j(z)=0,\quad \mbox{for all } k\geq0,
\end{equation}
then $\mu_j=0$ for all $j\geq1$.
\end{lemma}

\begin{proof}
Here we argue by induction. For $j=1$, by dividing \eqref{lem:eq:ann} by $r_1^k$, we have
\begin{equation}\label{lem:eq:ann0}
\int_{\mathbb T}z^k\,d\mu_1(z)+\sum_{j=2}^\infty\left(\frac{r_j}{r_1}\right)^k\int_{\mathbb T}z^k\,d\mu_j(z)=0.
\end{equation}
For $0<t<1$, denote by $P_t$ the Poisson kernel on $\mathbb T$. Recall that
\begin{equation*}
\int_{\mathbb T}z^k\,d(P_t*\mu_j)(z)=t^k\int_{\mathbb T}z^k\,d\mu_j(z),\quad k\geq0,\ j\geq1.
\end{equation*}

Define
\begin{equation*}
\sigma_1=\sum_{j=2}^\infty P_{r_j/r_1}*\mu_j,
\end{equation*}
then
\begin{equation*}
\int_{\mathbb T}z^k\,d(\mu_1+\sigma_1)(z)=0,\quad k\geq0.
\end{equation*}
Thus $\mu_1+\sigma_1\ll m$, the normalized Lebesgue measure on $\TT$, due to the F. and M. Riesz Theorem. Moreover, note that the Poisson convolution is contractive in total variation (for example see \cite[Theorem 6.4(a)]{ABR}), then
\begin{equation*}
\sum_{j=2}^\infty\|P_{r_j/r_1}*\mu_j\|\leq\sum_{j=2}^\infty\|\mu_j\|<\infty
\end{equation*}
which in turn implies that each term is absolutely continuous with respect to $m$. Hence $\sigma_1\ll m,$ and then $\mu_1\ll m.$ Therefore $\mu_1=0$ by the hypotheses that $\mu_1$ is purely atomic.

To be proceed, suppose $\mu_1=\cdots=\mu_{j-1}=0.$ Then \eqref{lem:eq:ann0} can be written as
\begin{equation*}
\int_{\mathbb T}z^k\,d\mu_{j}(z)+\sum_{\ell>j}\left(\frac{r_\ell}{r_j}\right)^k\int_{\mathbb T}z^k\,d\mu_\ell(z)=0.
\end{equation*}
By the same argument with
\begin{equation*}
\sigma_j=\sum_{\ell>j}P_{r_\ell/r_j}*\mu_\ell,
\end{equation*}
one can show that $\mu_j\ll m$, hence $\mu_j=0$.
\end{proof}

\begin{theorem}\label{ex-diag-mc}
Let $\lambda=(\lambda_1,\cdots,\lambda_d)$ such that $0\leq|\lambda_1|\leq\cdots\leq|\lambda_q|<|\lambda_{q+1}|=\cdots=|\lambda_d|=1,$ and $\lambda^\alpha\neq1$ $(\alpha\in\mathbb Z^d\setminus\{\boldsymbol 0\}).$
Denote $\mathbb N_0=\{0,1,2,\ldots\}$, and put $d_\alpha=\lambda^\alpha$ for the multi-index $\alpha\in\mathbb N_0^d$. Then the corresponding diagonal operator $D=\diag(d_\alpha)$ has the minimal commutant property on a separable Hilbert space.
\end{theorem}

\begin{proof}
Since $\lambda^\alpha\neq1$ $(\alpha\in\mathbb Z^d\setminus\{\boldsymbol 0\})$, then $d_\alpha=\lambda^\alpha$ are pairwise distinct  for multi-indices $\alpha\in\mathbb N_0^d$, and thus, by Theorem \ref{thm-diag}, it suffices to show that
\begin{equation*}
\left\{c\in\ell^1:\sum_{\alpha\in\mathbb N_0^d}c_\alpha d_\alpha^k=0\quad\forall k\geq0\right\}=\{0\}.
\end{equation*}

We first assume that $q\geq1$ and write $r_j=|\lambda_j|<1,$ $1\leq j\leq q.$  For $\alpha=(\beta,\gamma)$ with
$\beta\in\NN_0^q$ and $\gamma\in\NN_0^{d-q},$ we have
\begin{equation*}
|d_\alpha|=\prod_{j=1}^q r_j^{\beta_j}.
\end{equation*}
Consider the set of distinct radii
\begin{equation*}
\mathcal R=\left\{\prod_{j=1}^q r_j^{\beta_j}:\beta\in\mathbb N_0^q\right\}.
\end{equation*}
Since for every $\varepsilon>0$, there are only finitely many $\beta\in\mathbb N_0^q$ satisfying
\begin{equation*}
\prod_{j=1}^q r_j^{\beta_j}\geq\varepsilon.
\end{equation*}
Then the set $\mathcal R$ can be enumerated as
\begin{equation*}
R_0>R_1>R_2>\cdots>0,\quad R_j\rightarrow0.
\end{equation*}

Given $c=(c_\alpha)\in\ell^1$, for each $R\in\mathcal R$, define
\begin{equation*}
\mu_R=\sum_{\substack{\alpha\in\mathbb N_0^d\\
|d_\alpha|=R}}c_\alpha\delta_{d_\alpha/R},
\end{equation*}
which is a measure on $\TT$. 
Thus $\mu_R$ is purely atomic and
\begin{equation*}
\|\mu_R\|=\sum_{|d_\alpha|=R}|c_\alpha|.
\end{equation*}
Also we have
\begin{equation*}
\sum_{R\in\mathcal R}\|\mu_R\|\leq\sum_{\alpha\in\mathbb N_0^d}|c_\alpha|<\infty.
\end{equation*}

Moreover, we can rewrite
$$\begin{aligned}
0&=\sum_{\alpha\in\mathbb N_0^d}c_\alpha d_\alpha^k\\
&=\sum_{R\in\mathcal R}
R^k
\sum_{|d_\alpha|=R}
c_\alpha\left(\frac{d_\alpha}{R}\right)^k\\
&=\sum_{R\in\mathcal R}
R^k\int_{\mathbb T}z^k\,d\mu_R(z),
\qquad k\geq0.
\end{aligned}$$
So $\mu_R=0$ for all $R\in\mathcal R$ by Lemma \ref{lem:ann}. Since the atoms within each $\mu_R$ are distinct, then $c_\alpha=0$ for every $\alpha\in\mathbb N_0^d.$

Next we consider the case $q=0$. Then all $d_\alpha$ are pairwise distinct and lie on $\TT$. Similarly, define the purely atomic
\begin{equation*}
\mu=\sum_{\alpha\in\mathbb N_0^d}c_\alpha\delta_{d_\alpha}.
\end{equation*}
Given $c=(c_\alpha)\in\ell^1$, again we obtain
\begin{equation*}
\|\mu\|=\sum_\alpha|c_\alpha|<\infty.
\end{equation*}
In this case, if $\sum_{\alpha\in\mathbb N_0^d}c_\alpha d_\alpha^k=0$, then
\begin{equation*}
\int_{\mathbb T}z^k\,d\mu(z)=0,\quad k\geq0.
\end{equation*}
By the F. and M. Riesz theorem, $\mu\ll m$. Since $\mu$ is purely atomic, $\mu=0$, and hence $c_\alpha=0$ for every $\alpha$.
\end{proof}

Let $d\geq 1$, let us endow $\mathbb{C}^{d}$ with the standard inner product $\langle z,w\rangle=\sum_{j=1}^{d}z_{j} \overline{w_{j}}$ and let us denote $|z|=\sqrt{\langle z,{z}\rangle}$.
The classical Fock space on $\mathbb{C}^{d}$ is the separable Hilbert space of holomorphic functions on $\CC^d$ defined by
\begin{equation*}
\fcd=\left\{f\in H(\mathbb{C}^{d}):  \|f\|^{2}:=\frac{1}{(2\pi)^{d}}\int_{\mathbb{C}^{d}}|f(z)|^{2}e^{-{|z|^{2}}/{2}}dA(z)<\infty\right\}.
\end{equation*}
The set of polynomials $\{z^{\alpha}:\alpha\in \mathbb{Z}_+^{d}\}$ is an orthogonal basis of $\fcd$ with $\|z^{\alpha}\|^{2}=2^{|\alpha|}\alpha!$
for all $\alpha\in\ZZ_+^d.$
A holomorphic map $\varphi:\mathbb{C}^{d}\rightarrow\mathbb{C}^{d}$
induces a composition operator $C_\varphi$ defined by $C_{\varphi}f=f\circ\varphi$ for any  $f\in \fcd$. Carswell et al \cite{CMS} have characterized when $C_{\varphi}$ is a bounded or compact composition operator on $\fcd$, as follows.

\begin{theorem}\label{bound-com}
Let $\varphi:\mathbb{C}^{d}\rightarrow\mathbb{C}^{d}$ be holomorphic.
\begin{itemize}
  \item[(a)] $C_{\varphi}$ is bounded on $\fcd$ if and only if $\varphi(z)=Az+B$, where $A\in \mathcal M_d(\mathbb C)$, $\|A\|\leq1$, $B\in \mathbb{C}^{d}$,
and $\langle A\zeta,B\rangle=0$ for all $\zeta\in\mathbb{C}^{d}$ with $|A\zeta|=|\zeta|$.
  \item[(b)] $C_{\varphi}$ is compact on $\fcd$ if and only if $\varphi(z)=Az+B$, where $A\in \mathcal M_d(\mathbb C)$ and $\|A\|<1$, $B\in \mathbb{C}^{d}$.
\end{itemize}
\end{theorem}

Now we can characterize completely the minimal commutant property of composition operators induced by $\varphi(z)=Az$ with a normal matrix $A$, which also solves partially that problem arising from \cite[Remark 5.3]{BY}.

\begin{corollary}\label{MC-com-Fock}
Let $\varphi(z)=Az$ with a normal matrix $A$ and $\|A\|\leq1$. Then $C_{\varphi}$ has the minimal commutant property if and only if $(\lambda_j)$, the sequence of eigenvalues of $A$, satisfies that $\lambda^\alpha\neq1$ $(\alpha\in\mathbb Z^d\setminus\{\boldsymbol 0\})$.
\end{corollary}

\begin{proof}
Since $A$ is normal, then $C_{\varphi}$ is unitarily similar to the diagonal operator $\Lambda:$
$z^{\alpha}\mapsto\lambda^{\alpha}z^{\alpha}$ for all $\alpha\in \mathbb{N}^{d}$, and  $\sigma_{p}(C_{\varphi})=\sigma_{p}(\Lambda)=\left\{\lambda^{\alpha}:\alpha\in \mathbb{N}_0^{d}\right\}$.

The necessity is obtained immediately from \cite[Theorem 3.1]{BY}. Conversely, assume $\lambda^\alpha\neq1$ $(\alpha\in\mathbb Z^d\setminus\{\boldsymbol 0\})$.
If $\rho(A)<1$, i.e., $|\lambda_j|<1$ for all $j$, then $C_{\varphi}$ has the minimal commutant property by \cite[Theorem 5.2]{BY}. If $\rho(A)=1$, then $A$ has at least one eigenvalue of modulus $1$. In this case, again $C_{\varphi}$ has a minimal commutant by Theorem \ref{ex-diag-mc}.
\end{proof}

\medskip
\subsection{The case containing one Jordan block}~

\medskip
To consider the case containing one Jordan block of order $2$ inspired by \cite[Theorem 3.1]{BY}, we need the next two auxiliary lemmas.

\begin{lemma}\label{lem:ac-mear}
Let $\nu$ be a finite complex Borel measure on $\mathbb T$. If
\begin{equation*}
\big(\widehat\nu(n)\big)_{n\geq0}\in\ell^2,\quad
\widehat\nu(n)=\int_{\mathbb T}z^n\,d\nu(z),
\end{equation*}
then $\nu\ll m,$ the normalized Lebesgue measure on $\TT$.
\end{lemma}

\begin{proof}
Since the nonnegative Fourier coefficients of $\nu$ belong to $\ell^2$, there exists $f\in L^2(\mathbb T)$ such that
\begin{equation*}
\int_{\TT}z^n f(z)\,dm(z)=\widehat\nu(n),\qquad n\geq0,
\end{equation*}
and
\begin{equation*}
\int_{\TT}z^n f(z)\,dm(z)=0,\quad n<0.
\end{equation*}
So the measure $d\rho=(d\nu-f\,dm)$ satisfies
\begin{equation*}
\int_{\mathbb T}z^n\,d\rho(z)=0,\quad n\geq0.
\end{equation*}
By the F. and M. Riesz theorem, $\rho\ll m$. Hence $\nu\ll m$.
\end{proof}

\begin{lemma}\label{thm:radial-layer-elimination}
Let $\mathcal I$ be a countable index set. For each $\alpha\in\mathcal I$, let $m_\alpha\in\mathbb N_0$, $\mu_\alpha=r_\alpha\xi_\alpha$ with $0<r_\alpha\leq 1$ and $\xi_\alpha\in\mathbb T.$
Assume that the following conditions hold:
\begin{itemize}
\item[(i)] The set of radii $\mathcal R=\{r_\alpha:\alpha\in\mathcal I\}$ is locally finite in $(0,1]$. Moreover, for every $r\in\mathcal R$, the integers $m_\alpha$ are uniformly bounded over $\{\alpha:r_\alpha=r\}$. And for $\alpha\in\mathcal I_r:=\{\alpha\in\mathcal I:r_\alpha=r\}$, the points $\xi_\alpha$  are pairwise distinct.

\item[(ii)] For every $\alpha\in\mathcal I$, there is a nilpotent operator $N_\alpha$ satisfying $A_\alpha=\mu_\alpha \exp(N_\alpha)$ and $N_\alpha^{m_\alpha+1}=0$. Meanwhile, there exist $a\in(0,1)$ and $C>0$ such that, for any $\alpha\in\mathcal I$ and $k\in\NN_0$,
    \begin{equation*}
    \|A_\alpha^k\|\leq [C(1+k)]^{m_\alpha}\max\{r_\alpha^k,a^{km_\alpha}\}.
    \end{equation*}
\end{itemize}
Given a family of trace-class operators $R_\alpha$ such that $\sum_{\alpha\in\mathcal I}\|R_\alpha\|_1<\infty$, if
\begin{equation*}
\sum_{\alpha\in\mathcal I}\Tr(R_\alpha A_\alpha^k)=0,\ k\in\NN_0,
\end{equation*}
then $\Tr(R_\alpha N_\alpha^j)=0$ for $\alpha\in\mathcal I$ and $0\leq j\leq m_\alpha.$
\end{lemma}

\begin{proof}
By the nilpotency of $N_\alpha$, we have
\begin{equation*}
A_\alpha^k=\mu_\alpha^k e^{kN_\alpha}=\mu_\alpha^k\sum_{j=0}^{m_\alpha}\frac{k^j}{j!}N_\alpha^j.
\end{equation*}
which in turn implies that the hypotheses $\sum_{\alpha\in\mathcal I}\Tr(R_\alpha A_\alpha^k)=0$ becomes
\begin{equation}\label{eq:ann0}
 \sum_{\alpha\in\mathcal I}\mu_\alpha^k\sum_{j=0}^{m_\alpha}k^j b_{\alpha,j}=0,\quad k\geq0.
\end{equation}
Here, and in the sequel, we denote $b_{\alpha,j}=\frac{1}{j!}\Tr(R_\alpha N_\alpha^j),$ $0\leq j\leq m_\alpha.$

Fix $r\in\mathcal R$ and set $M_r=\max_{\alpha\in\mathcal I_r}m_\alpha<\infty$ by the hypothesis (i). For $0\leq j\leq M_r$, define the measure
\begin{equation*}
\nu_{r,j}=\sum_{\substack{\alpha\in\mathcal I_r\\ m_\alpha\geq j}}b_{\alpha,j}\delta_{\xi_\alpha}.
\end{equation*}
Since only finitely many coefficients $b_{\alpha,j}$ occur there by the nilpotency of $N_\alpha$, and $|b_{\alpha,j}|\leq\frac{\|N_{\alpha}^{\,j}\|}{j!}\|R_\alpha\|_1$, then $\nu_{r,j}$ is a finite purely atomic measure.

Suppose that all coefficients belonging to radial layers strictly larger than $r$ have already been shown to vanish. Set
\begin{equation*}
E_k(r)=\sum_{\substack{\alpha\in\mathcal I\\r_\alpha\leq r}}\Tr(R_\alpha A_\alpha^k).
\end{equation*}
Choose $L\geq1$ sufficiently large that $a^L<r.$ Split $E_k(r)$ into
\begin{equation*}
E_k(r)=E_k^{(1)}(r)+E_k^{(2)}(r),
\end{equation*}
where the first sum is over $m_\alpha<L$ and the second over $m_\alpha\geq L$. It follows from the hypothesi (i), there exists $s_0\in(0,r)$ such that $r_\alpha\leq s_0$ whenever $r_\alpha<r$ and $m_\alpha<L.$ Then, by the hypothesis (ii),
$$|E_k^{(1)}(r)|\leq\sum_{\substack{r_\alpha\leq r\\m_\alpha<L}}\|R_\alpha\|_1\|A_\alpha^k\|\leq
 [C(1+k)]^L s_0^k\sum_{\alpha\in\mathcal I}\|R_\alpha\|_1.$$
For the second part, choose $k$ sufficiently large so that $q_k:=C(1+k)a^k<1.$ Then, for $m_\alpha\geq L$, $\|A_\alpha^k\|\leq q_k^{m_\alpha}\leq q_k^L.$ Hence, again by the hypothesis (ii),
$$|E_k^{(2)}(r)|\leq\sum_{\substack{r_\alpha\leq r\\m_\alpha\geq L}}\|R_\alpha\|_1\|A_\alpha^k\|\leq [C(1+k)]^L a^{Lk}\sum_{\alpha\in\mathcal I}\|R_\alpha\|_1.$$
All in all, we get
\begin{equation*}
|E_k(r)|\leq K_0(1+k)^L(s_0^k+a^{Lk}).
\end{equation*}
Since $s_0<r$ and $a^L<r$, there exists
$\rho\in(0,1)$ such that $\max\{s_0,a^L\}<r\rho.$ Then
\begin{equation*}
\frac{E_k(r)}{r^k}=O((1+k)^L\rho^k),
\end{equation*}
and therefore
\begin{equation*}
\frac{E_k(r)}{r^k}=o(k^{-N})\quad\text{for every }N\geq0.
\end{equation*}

By the induction hypothesis, \eqref{eq:ann0} reduces to
\begin{equation*}
r^k\sum_{j=0}^{M_r}k^j\widehat{\nu}_{r,j}(k)+E_k(r)=0,
\end{equation*}
where $\widehat{\nu}_{r,j}(k)=\int_{\TT}z^k d\nu_{r,j}(z).$ Thus
\begin{equation}\label{eq:ann0-Mr-induction}
\sum_{j=0}^{M_r}k^j\widehat{\nu}_{r,j}(k)=-\frac{E_k(r)}{r^k}.
\end{equation}
Dividing the above equation by $k^{M_r}$, together with the boundedness of the Fourier coefficients of finite measures, shows that
\begin{equation*}
\widehat{\nu}_{r,M_r}(k)=O(k^{-1}),\quad k\geq1.
\end{equation*}
So $\big(\widehat{\nu}_{r,M_r}(k)\big)_{k\geq0}\in\ell^2.$ It follows that $\nu_{r,M_r}$ is absolutely continuous with respect to the normalized Lebesgue measure $m$ on $\TT$ by Lemma \ref{lem:ac-mear}. Noting that it is a pure atomic measure, then $\nu_{r,M_r}=0.$ Substituting this into \eqref{eq:ann0-Mr-induction} and repeating the same
argument with $M_r-1,M_r-2,\ldots,0$ gives
\begin{equation*}
\nu_{r,M_r-1}=\cdots=\nu_{r,0}=0.
\end{equation*}
Since the points $\xi_\alpha$, $\alpha\in\mathcal I_r$, are pairwise distinct, then
\begin{equation*}
b_{\alpha,j}=0,\quad \alpha\in\mathcal I_r,\quad0\leq j\leq m_\alpha.
\end{equation*}

The hypotheses (i) implies that the distinct radial values can be enumerated as $r_1>r_2>\cdots,$ $r_k\longrightarrow0$. Starting with the largest radial layer and by the induction argument,
\begin{equation*}
b_{\alpha,j}=0
 \quad\text{for all }\alpha\in\mathcal I,\quad
 0\leq j\leq m_\alpha.
\end{equation*}
In other wrods, $\Tr(R_\alpha N_\alpha^j)=0.$
\end{proof}

Finally, we can show the minimal commutant of composition operators of the other case containing one Jordan block of order $2$.

\begin{corollary}\label{MC-com-J-Fock}
Let $\phi(z)=Bz$ with $\|B\|\leq1$,
\begin{equation*}
    B=\begin{pmatrix} \lambda_1 & \lambda_1 & & & &  \\
   & \lambda_1 & & & & \\
    &  &  \lambda_2 & & & \\
    &  &   & \ddots & & \\
   & & & & & \lambda_{d-1} \end{pmatrix},
\end{equation*}
where $d>1$, $|\lambda_1|<1$, $|\lambda_j|\leq1$ and  $\lambda_{d-1}\in\TT$. If $\lambda_1^{\alpha_1}\cdots\lambda_{d-1}^{\alpha_{d-1}}\neq1$ for any multi-index $\alpha\in\mathbb Z^{d-1}\setminus\{\boldsymbol 0\}$, then $C_{\phi}$ has a minimal commutant on the Fock space.
\end{corollary}

\begin{proof}
For $n\geq0$, define
\begin{equation*}
\Ee_n=\s\{z_1^{n-j}z_2^j:0\leq j\leq n\}.
\end{equation*}
Write $\alpha=(n,\beta)\in\NN_0^{d-1}$ with $n\in\NN_0$ and $\beta=(\beta_2,\ldots,\beta_{d-1})\in\NN_0^{d-2},$ denote
$z'^\beta=z_3^{\beta_2}\cdots z_d^{\beta_{d-1}},$ $\eta_\beta=\prod_{j=2}^{d-1}\lambda_j^{\beta_j}$ and $\mu_\alpha=\lambda_1^n\eta_\beta$. Then by Fubini's theorem,
\begin{equation*}
\fcd=\bigoplus_{\alpha=(n,\beta)\in\NN_0^{d-1}}\mathcal H_\alpha\ \mbox{ with } \mathcal H_{\alpha}=\Ee_n z'^\beta
\end{equation*}
is an orthogonal decomposition. In particular every $\mathcal H_{\alpha}$ reduces $C_{\phi}$, and denote $A_\alpha=C_{\phi}|_{\mathcal H_{\alpha}}$. Then
\begin{equation}\label{eq:com-dem}
C_{\phi}=\bigoplus_{\alpha=(n,\beta)\in\NN_0^{d-1}} A_\alpha.
\end{equation}
Define $D_n=z_2\frac{\partial}{\partial z_1}\big|_{\Ee_n}$. Then it is easy to see that $D_n(z_1^{n-k}z_2^k)=(n-k)z_1^{n-k-1}z_2^{k+1}$ and $D_n^{n+1}=0$, and thus
\begin{equation*}
\exp(kD_n)f(z_1,z_2)=f(z_1+kz_2,z_2)\ \ \mbox{ for } f\in\Ee_n.
\end{equation*}
Now a direct computation gives
\begin{equation*}
A_\alpha=C_\phi|_{\Hh_{\alpha}}=\mu_\alpha \exp(D_n),
\end{equation*}
and for every integer $k\geq0$,
\begin{equation}\label{eq:A-alpha-k}
A_\alpha^k=\mu_\alpha^k \exp(kD_n)=\mu_\alpha^k\sum_{j=0}^{n}\frac{k^j}{j!}D_n^j.
\end{equation}
Moreover, following from \cite[Proposition 4.5]{BWY}, $\ker(C_{\phi}-\mu_\alpha I)^{n+1}=\Hh_{\alpha},$ and then $\{A_\alpha\}'=\alg\{D_n\}=\s\{D_n^{0},\cdots,D_n^{n}\}$. 

Since
$\sigma(A_{\alpha})=\left\{\lambda_1^n\eta_\beta\right\}$, and $\lambda_1^n\eta_\beta=\lambda_1^m\eta_\gamma$ if and only if $n=m$ and $\beta=\gamma$, then by the similar argument in Lemma \ref{sum-com-inv}, $T\in\{C_{\phi}\}'$ has the form
\begin{equation*}
  \bigoplus_{\alpha=(n,\beta)\in\NN_0^{d-1}} T_\alpha\ \mbox{ with } \ T_\alpha \in\{A_\alpha\}'=\alg\{D_n\} \mbox{ and } \sup_\alpha\|T_\alpha\|<\infty.
\end{equation*}

As is well-known, every WOT-continuous linear functional admits the form
$\Phi(X)=\sum_{k=1}^m \langle X x_k,y_k\rangle$
for finitely many vectors $x_k,y_k\in \fcd$. Equivalently, there exists a finite-rank operator $R=\sum_{k=1}^m y_k\otimes x_k$ such that
$\Phi(X)=\Tr(RX).$ Therefore, it suffices to prove that: for any functional $\Phi$ satisfying $\Phi(p(A))=0$ for every
polynomial $p$, $\Phi$ vanishes on all of $\Mcal=\{C_{\phi}\}'$. For this, denote by $P_{\alpha}$ the orthogonal projection onto $\mathcal H_{\alpha}$ and define $R_{\alpha}=P_{\alpha}RP_{\alpha}.$ Also recall that the set of trace-class operators is a two-sided ideal in $B(\fcd)$, and the trace-class pinching inequality gives
\begin{equation*}
\sum_\alpha\|R_\alpha\|_1\leq \|R\|_1.
\end{equation*}

To be proceed, we assume that $\Phi(C_\phi^k)=0$, i.e.,
\begin{equation*}
\Tr(R C_\phi^k)=\sum_{\alpha}\Tr(R_\alpha A_\alpha^{k})=0,\ \mbox{ for any } k\geq0
\end{equation*}
by \eqref{eq:com-dem}. Together \eqref{eq:A-alpha-k}, it follows that
\begin{equation}\label{eq:moment}
\sum_{\substack{\alpha\in\NN_0^{d-1}\\ \alpha=(n,\beta)}}\mu_\alpha^k\sum_{j=0}^{n}k^j b_{\alpha,j}=0,\quad k\geq0,
\end{equation}
where $b_{\alpha,j}=\frac{1}{j!}\Tr(R_{\alpha} D_{n}^{j})$ for $\alpha=(n,\beta)$ and $0\leq j\leq n$. Next we intend to use Lemma \ref{thm:radial-layer-elimination} to obtain $\Tr(R_{\alpha} D_{n}^{j})=0$.

To this end, let $J=\{j:|\lambda_j|<1\}\subseteq\{1,\ldots,d-1\},$ $U=\{1,\ldots,d-1\}\setminus J.$ For each multi-index $\alpha$, set
\begin{equation*}
r_\alpha=|\mu_\alpha|=\prod_{j\in J}|\lambda_j|^{\alpha_j}.
\end{equation*}
The value $r_\alpha$ is independent of the coordinates $\alpha_j$ with
$j\in U$. For every $\varepsilon>0$, there are only finitely many choices of
$(\alpha_j)_{j\in J}$ for which $r_\alpha\geq\varepsilon$. So the set
\begin{equation*}
\mathcal R=\{r_\alpha:\alpha\in\NN_0^{d-1}\}
\end{equation*}
is a countable subset of $(0,1]$, is locally finite in $(0,1]$, and
can be enumerated as $1=r^{(1)}>r^{(2)}>\cdots\rightarrow0.$ Given any $r\in\mathcal R$, define $M_r=\max\{\alpha_1:r_\alpha=r\}.$ Set $\xi_\alpha=\frac{\mu_\alpha}{r_\alpha}\in\TT$, and the points $\xi_\alpha$ are pairwise distinct by the hypotheses. 

Following from \cite[Theorem 6.1]{BWY26}, $C_\phi$ is power-bounded, and then so is $A_\alpha=C_\phi|_{\Hh_{\alpha}}$ for $\alpha=(n,\beta)$. Here a more precise asymptotic estimate is available. In fact, note that $\begin{pmatrix}1&1\\0&1\end{pmatrix}^k=\begin{pmatrix}1&k\\0&1\end{pmatrix}$, and thus there is a constant $C>0$ such that
\begin{equation*}
\|A_\alpha^k\|\leq\big(C(1+k)|\lambda_1|^k\big)^{n}|\eta_\beta|^{k}.
\end{equation*}
Therefore, we obtain all conditions of Lemma \ref{thm:radial-layer-elimination}, and then $\Tr(R_{\alpha} D_{n}^{j})=0$.


Finally, given any $T\in\{C_{\phi}\}'$, write
\begin{equation*}
T=\bigoplus_{\alpha=(n,\beta)\in\NN_0^{d-1}} T_\alpha \ \mbox{ with } T_\alpha=\sum_{j=0}^{n}t_{\alpha,j}D_{n}^{j}.
\end{equation*}
Note that
\begin{equation*}
\Tr(R_\alpha T_\alpha)=\sum_{j=0}^{\alpha_1}t_{\alpha,j}\Tr(R_\alpha D_{n}^{j})=0,
\end{equation*}
and
\begin{equation*}
|\Phi(T)|\leq\sum_\alpha|\Tr(R_\alpha T_\alpha)|\leq\|T\|\sum_\alpha\|R_\alpha\|_1<\infty.
\end{equation*}
Therefore $\Phi(T)=0$ for every $T\in\Mcal$.
\end{proof}

\medskip
\noindent \textbf{Acknowledgment.}
M. Wang and X. Yao were supported by National Science Foundation of China (No. 12571145). The authors would like to thank Professor Chunlan Jiang for his insightful and helpful comments.

\end{document}